\documentclass[11pt]{article}

\usepackage[a4paper,margin=1in]{geometry}
\usepackage[T1]{fontenc}
\usepackage{lmodern}
\usepackage{microtype}
\usepackage{amsmath,amssymb,amsthm,amscd,mathtools}
\mathtoolsset{showonlyrefs}
\usepackage{mathrsfs}
\usepackage{enumitem}
\usepackage{xcolor}
\usepackage{xspace}
\usepackage[colorlinks=true,linkcolor=blue!55!black,citecolor=blue!55!black,urlcolor=blue!55!black]{hyperref}
\usepackage[backend=bibtex,style=alphabetic,maxnames=10,minnames=10,maxbibnames=99,minalphanames=4,maxalphanames=10]{biblatex}
\newtheorem{theorem}{Theorem}[section]
\newtheorem{maintheorem}{Theorem}

\newtheorem{proposition}[theorem]{Proposition}
\newtheorem{lemma}[theorem]{Lemma}
\newtheorem{corollary}[theorem]{Corollary}

\theoremstyle{definition}
\newtheorem{definition}[theorem]{Definition}
\newtheorem{remark}[theorem]{Remark}

\newcommand{\remarkendmark}{\ensuremath{\diamond}}
\let\proofendmark\qedsymbol
\AddToHook{env/remark/begin}{%
  \renewcommand{\qedsymbol}{\remarkendmark}\pushQED{\qed}}
\AddToHook{env/remark/end}{\popQED}
\AddToHook{env/proof/begin}{\let\qedsymbol\proofendmark}

\newtheorem{example}[theorem]{Example}

\newcommand{\Qp}{\mathbb Q_p}
\newcommand{\Zp}{\mathbb Z_p}

\newcommand{\cO}{\mathcal O}

\newcommand{\cM}{\mathcal M}
\newcommand{\cF}{\mathcal F}
\newcommand{\cL}{\mathcal L}

\newcommand{\Pone}{\mathbb P^1}
\newcommand{\Period}{\mathbb P_C^{n-1}}

\newcommand{\la}{\mathrm{la}}

\newcommand{\an}{\mathrm{an}}
\newcommand{\proet}{\mathrm{pro\acute et}}
\newcommand{\et}{\mathrm{\acute et}}

\newcommand{\Ind}{\operatorname{Ind}}

\newcommand{\JL}{\operatorname{JL}}
\newcommand{\Lie}{\operatorname{Lie}}

\newcommand{\Rep}{\operatorname{Rep}}

\newcommand{\acong}{\mathrel{\overset{a}{\cong}}}
\newcommand{\ain}{\mathrel{\overset{a}{\in}}}
\DeclareMathOperator{\coker}{coker}

\title{On the cohomological amplitude of the locally analytic Jacquet--Langlands functor}
\author{Yu Xiao}
\date{\today}

\begin{document}
\maketitle

\begin{abstract}
We prove a local cohomological bound
\(R\Gamma_{\proet}(X,\mathcal F)\in D^{[0,\dim X]}\)
for a pro-\'etale sheaf with a derived \(p\)-complete integral
lattice whose modulo $p^m$ reductions are filtered colimits of finite local
systems on an affinoid topologically of finite type over $C$. 
For \(G=\mathrm{GL}_n(F)\), with \(n\geq2\) and \(F/\Qp\)
finite, we apply these local bounds to variants of Scholze's Jacquet--Langlands functor  
$\JL_\mathscr{B}(\Pi)\coloneqq R\Gamma_{\proet}(\mathbb{P}^{n-1}_C,\mathcal{F}_{\Pi}\otimes^{L, \square}_{\mathbb{Q}_p}\mathscr{B})$
for $\mathscr{B}\in\{\mathbb{Q}_p, \widehat{\mathscr{O}}, \mathbb{B}_J\}$ 
and Banach or compact-type LB representations $\Pi$ of $G$, obtaining
\(\JL_{\mathscr{B}}(\Pi)\in D^{[0,2(n-1)]}\) for each choice of $\mathscr{B}$. 
We also prove that for arbitrary continuous principal series $I(\eta)=(\mathrm{Ind}_B^G\eta)^{\mathrm{cont}}$, 
we have $\JL(\Ind_B^G\eta)\in D^{[0, n]}$, while $\JL_{\widehat{\mathscr{O}}}(I(\eta)),\JL_{\mathbb{B}_J}(I(\eta))\in D^{[0, n-1]}$.
\end{abstract}

\begingroup
\hypersetup{linkcolor=black}
\tableofcontents
\endgroup
\clearpage

\section*{Introduction}
\addcontentsline{toc}{section}{Introduction}

Let \(F/\Qp\) and \(E/\Qp\) be finite extensions and let
\(n\geq2\). Put
\begin{equation*}
 C=\widehat{\overline F},\qquad
 G=\mathrm{GL}_n(F),\qquad
 X=\Period.
\end{equation*}
Let \(D/F\) be the central division algebra of invariant \(1/n\)
and let $B\subset G$ be the upper-triangular Borel subgroup.
The infinite-level Lubin--Tate space \(\cM\) has a
\(D^\times\)-equivariant pro-\'etale \(G\)-torsor
\begin{equation*}
 \pi=\pi_{\mathrm{GH}}:\cM\longrightarrow X,\qquad
 \cM/G\simeq X,
\end{equation*}
with commuting \(G\)- and \(D^\times\)-actions. 

Scholze \cite{Scholze} constructs a Jacquet--Langlands functor
for admissible $p$-power-torsion representations of $G$.
Passing to limits extends it to unitary admissible Banach
representations.  Dospinescu and Rodr\'iguez Camargo \cite{DRC}
use solid formalism to study locally analytic representations,
prove compatibility with taking locally analytic vectors, and
establish compatibility of infinitesimal characters with structural
coefficient sheaf.

In this article we extend the cohomological amplitude bounds
to Banach representations and LB representations of compact
type, without admissibility or unitarity assumptions.  Following
\cite{DRC}, for such a representation $\Pi$ we define a coefficient
sheaf on $X_{\proet}$ by
\[\cF_\Pi(W)\coloneqq C_G(|W\times_{\mathbb{P}^{n-1}}\cM|, \Pi), \qquad W\in\mathbb{P}^{n-1}_{\proet}.\]
We set \[\JL(\Pi)=R\Gamma_{\proet}(X, \cF_\Pi), \qquad \JL^i(\Pi)=H^i\JL(\Pi).\]

Similarly, for the structure sheaf $\widehat{\mathscr O}$ and period sheaves $\mathbb{B}_J$ used in \cite{DRC}, we define 
\[
\JL_{\widehat{\mathscr O}}(\Pi)=R\Gamma_{\proet}(X, \cF_\Pi\otimes_{\Qp}^{L,\square}\widehat{\mathscr O}),\quad
\JL_{\mathbb{B}_J}(\Pi)=R\Gamma_{\proet}(X, \cF_\Pi\otimes_{\Qp}^{L,\square}\mathbb{B}_J).
\]
Here the solid tensor product is equivalent to the $p$-adic completed tensor product for Banach representation and the colimit of $p$-adic completed tensor product for LB representations. 

We first prove the following local cohomological bound for a pro-\'etale sheaf with a derived \(p\)-complete integral lattice whose reductions are filtered colimits of finite local systems on an affinoid.
\begin{theorem}[Theorem~\ref{thm:general-unitary-coefficient-amplitude}]
      Let $\mathrm{Spa}(A, A^+)$ be a connected affinoid, topologically of finite type over $C$, and let $\cF$ be a sheaf of $\Qp$-vector spaces on $\mathrm{Spa}(A, A^+)_{\proet}$, with an integral lattice $\cF^+\subset \cF$ satisfying
      \[\cF^+=R\varprojlim\cF^+/p^s,\qquad \cF^+[1/p]=\cF,\qquad \cF^+/p^s=v^*\cL_s\]
      for $\cL_s=\varinjlim_{\alpha}\cL_{s,\alpha}$ a filtered colimit of finite \'etale local systems on $\mathrm{Spa}(A, A^+)_{\et}$. Then we have
      \[R\Gamma_{\proet}(\mathrm{Spa}(A, A^+), \cF\otimes_{\Qp}^{L,\square}\mathscr{B})\in D^{[0, \dim A]}\] for all $\mathscr{B}\in\{\mathbb{Q}_p, \widehat{\mathscr O}, \mathbb B_J\}$.
\end{theorem}
Such a sheaf $\cF$ is called a complete unitary pro-\'etale local system. 
Applying this bound to the sheaves $\cF_{\Pi}$ and using a \v Cech complex, we get the following uniform cohomological bound:
\begin{theorem}[Theorems~\ref{thm:general-Banach-amplitude} and~\ref{thm:LB-amplitude}]
   The complexes above have cohomological amplitude contained in $[0, 2(n-1)]$. 
  \[\JL(\Pi), \JL_{\widehat{\mathscr O}}(\Pi), \JL_{\mathbb{B}_J}(\Pi)\in D^{[0, 2(n-1)]}(\mathrm{Rep}_E^{\square}(D^\times))\]
  for Banach or LB representations $\Pi$ of $G$. In particular, this bound holds for all admissible locally analytic representations of $G$. 
\end{theorem}
Here the cohomology amplitude $[a, b]$ of a complex $C^\bullet$ means $H^i(C^\bullet)=0$ for $i\notin[a, b]$. 
The bound is optimal already for unitary admissible Banach representations.

The work of Ludwig \cite{Ludwig} and Johansson--Ludwig--Hansen
\cite{JLH} gives sharper bounds for principal series with
torsion or unitary coefficients.  We extend these calculations
to general continuous inducing characters:

\begin{theorem}[Theorems~\ref{thm:continuous-principal-series-amplitude}
and~\ref{thm:principal-series-geometric-amplitude}]
  Let $I(\eta)=\Ind_B^G(\eta)$ be the continuous induction of a continuous character $\eta\colon B\to E^\times$. Then we have 
  \[\JL(I(\eta))\in D^{[0, n]}(\mathrm{Rep}_E^{\square}(D^\times)), \quad
  \JL_{\widehat{\mathscr O}}(I(\eta)),
  \JL_{\mathbb B_J}(I(\eta))
  \in D^{[0,n-1]}(\mathrm{Rep}_E^\square(D^\times)).
  \]
\end{theorem}
Compared with Ludwig's bound $[0, n-1]$, the extra degree $n$ comes from the unknown $R^1\lambda_*$ for a non-unitary pro-\'etale sheaf. 
For \(n=2\), the bound does not improve the general bound,
whereas the structural and period bounds improve to \([0,1]\).

We also prove the following compatibilities for $\mathscr{B}\in\{\mathbb{Q}_p, \widehat{\mathscr O}, \mathbb{B}_J\}$, with locally analytic vectors under
the acyclicity hypothesis \(\Pi^{RG\text{-}\la}\simeq\Pi^{G\text{-}\la}[0]\), which holds for admissible
Banach representations and admissible locally analytic representations.
\begin{itemize}
      \item{With locally analytic vectors:}  $\JL_{\mathscr{B}}(\Pi)^{RD^\times-\mathrm{la}}=\JL_{\mathscr{B}}(\Pi^{G-\mathrm{la}})$.
      \item{With twisting characters:} $\JL_{\mathscr{B}}(\Pi\otimes(\chi\circ\det))=\JL_{\mathscr{B}}(\Pi)\otimes((\chi\circ\mathrm{Nrd})\boxtimes(\chi\circ\mathrm{rec}_F)^{-1})$.
      \item{With central characters:} If $\Pi$ admits a central character $\zeta$, then so does $\JL_{\mathscr{B}}(\Pi)$ under $Z(G)=F^\times=Z(D^\times)$.
      \item{With infinitesimal characters:} If $\Pi^{G-\la}$ admits an infinitesimal character, then so does $\JL_{\widehat{\mathscr{O}}}(\Pi^{G-\la})$ via $Z(\Lie(D^\times)\otimes_{\mathbb{Q}_p}E)\cong Z(\Lie(G)\otimes_{\mathbb{Q}_p}E)$.
\end{itemize}

For general preservation of infinitesimal characters, we need more information of the primitive comparison map $\JL(\Pi)\otimes^{L,\square}_{\mathbb{Q}_p}C\to\JL_{\widehat{\mathscr{O}}}(\Pi)$.
But in general this map is not a quasi-isomorphism, see Example~\ref{ex:standard-representation-primitive-failure} and compare with \cite[Example~1.1.2]{ALM}.

We use the local cohomology bound of Theorem~\ref{thm:general-unitary-coefficient-amplitude-A}
on each small affinoid chart. 
And then we construct a rational cover \(\{U_i\}_{i=0}^{n-1}\) of $\mathbb{P}^{n-1}_C$ on which $\cF_\Pi|_{U_i}$ satisfies the conditions in Theorem~\ref{thm:general-unitary-coefficient-amplitude-A}.
The spectral sequence
\[
 E_1^{a,b}=\bigoplus_{0\leq i_0<\cdots<i_a<n}
 H^b_{\proet}(U_{\{i_0,\ldots,i_a\}},\cF_\Pi\otimes_{\mathbb{Q}_p}^{L, \square}\mathscr{B})
 \Longrightarrow H^{a+b}_{\proet}(\mathbb P_C^{n-1},\cF_\Pi\otimes_{\mathbb{Q}_p}^{L, \square}\mathscr{B})=\JL_{\mathscr{B}}(\Pi)
\]
then yields the global interval \([0,2(n-1)]\).
We obtain the same local interval for structural and period
coefficients using the affinoid perfectoid universal pro-finite-\'etale cover
and almost purity for perfectoid spaces. We also use the fact that the
\'etale fundamental group has \(p\)-cohomological dimension at most \(\dim X\) using the \'etale $K(\pi, 1)$ comparison..

\paragraph{Outline of paper}In section~\ref{sec:general-local-amplitude}, we introduce the complete unitary pro-\'etale $\mathbb{Q}_p$-local system, and prove an Artin-Grothendieck type vanishing for such sheaf
together with the cohomologies after tensoring with $\widehat{\mathscr{O}}$ and $\mathbb{B}_J$ on affinoids over $C$ (Theorem~\ref{thm:general-unitary-coefficient-amplitude}). 

In section~\ref{sec:cohomological-amplitude}, we apply the local vanishing to the cover in Proposition~\ref{prop:GH-rational-cover} and get the corresponding Cohomological amplitude (Theorems~\ref{thm:period-amplitude} and \ref{thm:LB-amplitude}).

In section~\ref{sec:continuous-parabolic-induction}, we generalize the methods in \cite{Ludwig} and \cite{JLH} to non-unitary continuous principal series (Theorems~\ref{thm:continuous-principal-series-amplitude} and \ref{thm:principal-series-geometric-amplitude}).

\paragraph{Acknowledgement}The author is deeply grateful to his advisor, Professor Yongquan Hu, for suggesting the area of $p$-adic Jacquet-Langlands functor studied in this paper, 
and for his constant guidance and encouragement. 
The author would also like to thank Gabriel Dospinescu for helpful discussions. 
He pointed out the importance of the local cohomological bound for general Banach sheaves, 
and suggested the use of Bhatt-Mathew's torsion vanishing for \'etale local systems on smooth affinoids. 

\paragraph{AI disclosure}This article was drafted with the assistance of ChatGPT, using the model GPT-6 Astra, which helped with proofreading, correcting grammar and finding references.  
The author first wrote a complete draft with a weaker local bound from naive Artin vanishing $[0, 2\dim]$ for \'etale torsion local systems, 
and only deduced the global bound $[0,2n-1]$ for $\JL$ from the period fundamental sequence. What GPT did is to find the reference \cite{GabberZavyalov} about the Artin-Grothendieck vanishing result.
The author takes the full responsibility of this article. 

\section{Cohomological amplitude for complete unitary pro-\'etale sheaves on affinoids}
\label{sec:general-local-amplitude}

Let \(C/\Qp\) be a complete algebraically closed nonarchimedean
field. 
\begin{definition}
\label{def:complete-unitary-proetale-local-system}
Let $S$ be a rigid analytic variety over $C$ and let $\cF$ be a sheaf of $\mathbb{Q}_p$ vector spaces on the pro-\'etale site $S_{\proet}$.
We say that $\cF$ is a complete unitary pro-\'etale $\mathbb{Q}_p$-local system if there is a subsheaf $\cF^+\subset\cF$ such that
\[
 \mathcal F=\mathcal F^+[1/p],\qquad
 \mathcal F^+\xrightarrow{\sim}
 R\!\varprojlim_s
 \bigl(\mathcal F^+\otimes_{\Zp}^L\Zp/p^s\bigr).
\]and each $\cF^+/p^s$ is a filtered colimit of \'etale local systems: 
\[
 \mathcal F^+/p^s\mathcal F^+\cong\nu^*\mathcal L_s,
 \qquad
 \mathcal L_s=\varinjlim_\alpha\mathcal L_{s,\alpha},
\]where every \(\mathcal L_{s,\alpha}\) is a finite \'etale
\(\Zp/p^s\)-local system.
The subsheaf $\cF^+$ is called a lattice in $\cF$.
\end{definition}
We focus on affinoids over $C$ and use the following notation:
\[
 S=\operatorname{Spa}(A,A^\circ),\qquad d=\dim A,\qquad
 \nu:S_{\proet}\longrightarrow S_{\et},\qquad
 c_A:S_{\et}\longrightarrow(\operatorname{Spec}A)_{\et},
\]
where \(A\) is a \(C\)-affinoid algebra and we assume $S$ is connected. 
We first introduce the main result of this section:
\begin{maintheorem}
\label{thm:general-unitary-coefficient-amplitude-A}
Let \(\mathcal F\) be a complete unitary pro-\'etale \(\Qp\)-local system
on \(S\), with a lattice \(\mathcal F^+\) as in
Definition~\ref{def:complete-unitary-proetale-local-system}.
Then for every compact rational interval \(J\subset(0,\infty)\)
and every \(\mathscr B\in\{\Qp,\widehat{\mathscr O},\mathbb B_J\}\),
we have
\begin{equation*}
 R\Gamma_{\proet}(S,\mathcal F^+),\ 
 R\Gamma_{\proet}(S,\mathcal F\widehat{\otimes}_{\Qp}\mathscr B)
 \in D^{[0,d]}.
\end{equation*}
The completed tensor products agree with
\(\mathcal F\otimes_{\Qp}^{L,\square}\mathscr B\), and the
rational bounds hold in the solid category.
\end{maintheorem}
We restate and prove the separate cases in Theorem~\ref{thm:general-local-lattice-amplitude} and
Theorem~\ref{thm:general-unitary-coefficient-amplitude}. 

\subsection{Torsion vanishing on affinoids}
We begin with the vanishing result for \'etale torsion sheaves over affinoids,
and then relate it to the \'etale fundamental group by the $K(\pi, 1)$-realization \cite[Theorem~4.9]{ScholzeRigid}. 

\subsubsection{An Artin-Grothendieck vanishing}
In this subsection, we refer to \cite{GabberZavyalov} for algebraic sheaves and the vanishing result. 
An \'etale torsion sheaf \(\mathcal L\) is algebraic if
\(\mathcal L\cong c_A^*\mathcal G\) for a torsion sheaf
\(\mathcal G\) on \((\operatorname{Spec}A)_{\et}\).

\begin{proposition}[Gabber--Zavyalov]
\label{prop:algebraic-Artin}
For every algebraic torsion sheaf \(\mathcal L\) on \(S_{\et}\),
\begin{equation}\label{eq:general-algebraic-Artin}
 R\Gamma_{\et}(S,\mathcal L)\in D^{[0,d]}.
\end{equation}
The same conclusion holds for filtered colimits of algebraic
torsion sheaves, in particular for filtered colimits of finite
locally constant \(\Zp/p^s\)-sheaves.  
\end{proposition}
\begin{proof}
The theorems 
\cite[Theorems~3.1 and~4.4]{GabberZavyalov} prove the bound for algebraic torsion sheaves $\cL$. 
By \cite[Lemma~4.3(1)]{GabberZavyalov}, finite \'etale local systems are algebraic.  Since the \'etale site \(S_{\et}\) is coherent, we can take filtered colimit to get
\[
 H^q_{\et}\!\left(S,\varinjlim_\alpha\mathcal L_\alpha\right)
 \cong\varinjlim_\alpha H^q_{\et}(S,\mathcal L_\alpha)=0
 \qquad(q>d).
\]
\end{proof}

\begin{remark}
\label{rem:smooth-local-system-Artin}
When \(A\) is smooth, \cite[Remark~7.4(2)]{BhattMathewArc}
gives an alternative proof of Proposition~\ref{prop:algebraic-Artin}
when $\cL$ is a finite \'etale \(\Zp/p^s\)-local system on \(S_{\et}\).  
We follow the method in \cite[proof of Theorem~7.3]{BhattMathewArc} 
and \cite[Tags~0A3R and~03SH]{Stacks}.

Every finite \'etale \(Y\to S\) is a smooth affinoid of dimension
at most \(d\).  By \cite[Remark~7.4(2)]{BhattMathewArc} and
algebraic--analytic comparison \cite[Corollary~3.2.3]{HuberEtale},
\[
 H^q_{\et}(Y,\mathbf F_p)=0\qquad(q>d).
\]
For \(s=1\), choose a finite
\'etale Galois cover \(\widetilde{S}\to S\) trivializing \(\mathcal L\),
with Galois group \(G\), and a Sylow \(p\)-subgroup \(P\subset G\).  Put
\[
 f:\widetilde{S}/P\longrightarrow S,\qquad
 r=\deg f=[G:P],\qquad p\nmid r.
\]
Since every simple \(\mathbf F_p[P]\)-module is trivial of dimension
one, the local system \(f^*\mathcal L\) admits a finite filtration
\[
 0=\mathcal L_0\subset\mathcal L_1\subset\cdots
 \subset\mathcal L_m=f^*\mathcal L,\qquad
 \mathcal L_j/\mathcal L_{j-1}\cong\mathbf F_p.
\]
Thus \(H^q_{\et}(\widetilde{S}/P,f^*\mathcal L)=0\) for \(q>d\).
On the other hand, the restriction map and the trace map satisfy
\[
 H^q_{\et}(S,\mathcal L)\xrightarrow{f^*}
 H^q_{\et}(\widetilde{S}/P,f^*\mathcal L)\xrightarrow{\operatorname{Tr}_f}
 H^q_{\et}(S,\mathcal L),\qquad
 \operatorname{Tr}_f\circ f^*=r\operatorname{id}.
\]
Since \(r\in\mathbf F_p^\times\), restriction is split injective, hence $H^q_{\et}(S,\cL)=0$ for $q>d$.
This proves Proposition~\ref{prop:algebraic-Artin} for \(s=1\) and finite \'etale local systems $\cL$. 

For general \(s\), apply the result to the finite local systems
\(p^j\mathcal L/p^{j+1}\mathcal L\), \(0\leq j<s\).
The \'etale site $S_{\et}$ is coherent. Hence taking filtered colimit we can get
\[
 H^q_{\et}\!\left(S,\varinjlim_\alpha\mathcal L_\alpha\right)
 \cong\varinjlim_\alpha H^q_{\et}(S,\mathcal L_\alpha)=0
 \qquad(q>d).
\qedhere
\]
\end{remark}

\subsubsection{The \'etale fundamental group}
\label{subsec:etale-fundamental-group}
Later we will use the \'etale fundamental group to simplify the computations in different settings.
We first recall the definition and the $K(\pi, 1)$ comparison used to bound the cohomological dimension. 

By assumption \(S\) is connected, we fix a geometric point \(\bar x\) throughout,
and put \(\Gamma=\pi_1^{\et}(S,\bar x)\).  For a finite discrete
\(p\)-primary continuous \(\Gamma\)-module \(N\), set
\[
 H_N=\ker\bigl(\Gamma\longrightarrow\operatorname{Aut}(N)\bigr),
 \qquad \Delta_N=\Gamma/H_N,
 \qquad r_N:Y_N\longrightarrow S,
\]
where \(r_N\) is the pointed finite \'etale right
\(\Delta_N\)-torsor corresponding to \(H_N\).
Write \(\underline N_{Y_N}\) for the constant \'etale sheaf on \(Y_N\).
Define the finite local system associated to \(N\) by
\begin{equation}\label{eq:fundamental-group-local-system}
 \mathcal L_N=(r_{N,*}\underline N_{Y_N})^{\Delta_N},
 \qquad (g\cdot f)(y)=g\,f(yg).
\end{equation}
Equivalently, for \(U\to S\) \'etale,
\[
 \mathcal L_N(U)=
 \bigl\{f:|U\times_SY_N|\longrightarrow N\text{ locally constant}:
          f(yg)=g^{-1}f(y)\text{ for }g\in\Delta_N\bigr\}.
\]
Evaluation at the chosen lift of \(\bar x\) identifies
\((\mathcal L_N)_{\bar x}\cong N\) with its given \(\Gamma\)-action.

\begin{proposition}
\label{prop:affinoid-fundamental-group-dimension}
One has
\begin{equation}\label{eq:general-affinoid-fundamental-group-dimension}
 \operatorname{cd}_p\Gamma\leq d,
 \qquad\text{i.e. }
 H^q_{\mathrm{cont}}(\Gamma,N)=0
 \quad(q>d)
\end{equation}
for every discrete \(p\)-primary continuous \(\Gamma\)-module \(N\).
\end{proposition}
\begin{proof}
For a finite $\Gamma$ module \(N\), the \(K(\pi,1)\) comparison
\cite[Theorem~4.9]{ScholzeRigid} and
Proposition~\ref{prop:algebraic-Artin} immediately give us
\[
 H^q_{\mathrm{cont}}(\Gamma,N)
 \cong H^q_{\et}(S,\mathcal L_N)=0\qquad(q>d).
\]
For a discrete \(p\)-primary \(N\), every orbit \(\Gamma n\)
is finite since \(\Gamma\) is compact. Thus $N$ is a filtered colimit of finite $\Gamma$ submodules
\[
 N=\varinjlim_{\substack{S\subset N\\S\text{ finite}}}N_S,
 \qquad
 N_S=\sum_{n\in S}\Zp[\Gamma]n,
 \qquad |N_S|<\infty.
\]
Every continuous map \(\Gamma^q\to N\) has finite image, hence
\[
 C^\bullet_{\mathrm{cont}}(\Gamma,N)
 =\varinjlim_S C^\bullet_{\mathrm{cont}}(\Gamma,N_S).
\]
Since taking filtered colimits is exact, we get
\[
 H^q_{\mathrm{cont}}(\Gamma,N)
 \cong\varinjlim_S H^q_{\mathrm{cont}}(\Gamma,N_S)=0
 \qquad(q>d).
\]
Thus by definition \(\operatorname{cd}_p\Gamma\leq d\).
If \(S\) is smooth, we can instead use Remark~\ref{rem:smooth-local-system-Artin}.\end{proof}

\subsection{Completion and the local bound}
In this subsection, we use the torsion results and take inverse limits to obtain
a vanishing result for the integral subsheaf. We first prove a general lemma on inverse limits.
\begin{lemma}
\label{lem:general-completion-amplitude}
Let \(R\) be a commutative ring and \((K_s)_{s\geq1}\) an inverse system
in \(D(R)\), with distinguished triangles
\begin{equation}\label{eq:general-reduction-triangle}
 K_1\longrightarrow K_{s+1}\longrightarrow K_s
 \longrightarrow K_1[1]
 \qquad(s\geq1).
\end{equation}
If \(K_1\in D^{[0,e]}\) for some \(e\geq0\), then
\[
 K_s\in D^{[0,e]}\quad(s\geq1),\qquad
 R\!\varprojlim_s K_s\in D^{[0,e]}.
\]
\end{lemma}

\begin{proof}
By induction on $s$ using \eqref{eq:general-reduction-triangle}, we get \(K_s\in D^{[0,e]}\) for all $s$.
Since the index set is countable, we have $R\varprojlim_s K_s\in D^{[0, e+1]}$. What's more at the top degree, the last terms of the long exact sequence
\[
 H^e(K_{s+1})\longrightarrow H^e(K_s)
 \longrightarrow H^{e+1}(K_1)=0
\]
gives surjective transition maps at the top-degree.  The Milnor sequences
\cite[Tag~0D6K]{Stacks}
\begin{equation}\label{eq:general-lattice-Milnor}
 0\longrightarrow R^1\!\varprojlim_s H^{j-1}(K_s)
 \longrightarrow H^j\!\left(R\!\varprojlim_s K_s\right)
 \longrightarrow\varprojlim_s H^j(K_s)\longrightarrow0
\end{equation}
then prove the assertion, since
\(R^1\!\varprojlim_s H^e(K_s)=0\).
\end{proof}

Now we can prove the main theorem when $\mathscr{B}=\mathbb{Q}_p$:
\begin{theorem}
\label{thm:general-local-lattice-amplitude}
Let \(\mathcal F\) be a complete unitary pro-\'etale local system on $S_{\proet}$, with integral lattice $\cF^+$. 
Then
\begin{align}
 R\Gamma_{\proet}(S,\mathcal F^+)
 &\simeq R\!\varprojlim_s R\Gamma_{\proet}(S,\cF^+/p^s)
   \in D^{[0,d]}(\Zp),
 \label{eq:general-integral-local-amplitude}\\
 R\Gamma_{\proet}(S,\mathcal F)
 &\simeq\bigl(R\Gamma_{\proet}(S,\mathcal F^+)\bigr)[1/p]
   \in D^{[0,d]}(\Qp).
 \label{eq:general-rational-local-amplitude}
\end{align}
\end{theorem}

\begin{proof}
By assumption $\cF^+/p^s=\nu^*\cL_s$, and we denote \(K_s=R\Gamma_{\et}(S,\mathcal L_s)\).
Since $\cF^+$ is derived $p$-complete, we can use the torsion comparison
\cite[Corollary~3.17]{ScholzeRigid} to obtain 
\[
 R\Gamma_{\proet}(S,\mathcal F^+)
 \simeq R\!\varprojlim_s K_s.
\]
Since the sheaf \(\mathcal F^+\) is $p$-torsion free, the following sequence of d\'evissage 
\begin{equation}\label{eq:general-lattice-successive-reduction}
 0\longrightarrow\mathcal F^+/p\mathcal F^+
 \xrightarrow{p^s}\mathcal F^+/p^{s+1}\mathcal F^+
 \longrightarrow\mathcal F^+/p^s\mathcal F^+
 \longrightarrow0,
\end{equation}
is exact, hence the triangles \eqref{eq:general-reduction-triangle}.
By Proposition~\ref{prop:algebraic-Artin}, we have
\(K_1\in D^{[0,d]}\). We then apply
Lemma~\ref{lem:general-completion-amplitude}.
Since $S$ is an affinoid, it is quasi-compact and 
$S_{\proet}$ is a coherent site by \cite[Proposition~3.12]{ScholzeRigid},
\begin{equation}\label{eq:general-local-rationalization}
 R\Gamma_{\proet}(S,\mathcal F^+[1/p])
 \simeq R\Gamma_{\proet}(S,\mathcal F^+)[1/p].
\end{equation}
This deduces the bound since localization is exact.
\end{proof}

In order to prove the vanishing results for $\mathcal{F}\widehat{\otimes}_{\mathbb{Q}_p}\widehat{\mathscr{O}}$
and $\cF\widehat{\otimes}_{\mathbb{Q}_p}\mathbb{B}_J$, we need an almost version of lemma~\ref{lem:general-completion-amplitude}. 
We use the following conventions. For an almost setting \((R,\mathfrak a)\), let \((-)^a\) denote
localization by the modules \(M\) with \(\mathfrak aM=0\).  We write
\[
 \begin{aligned}
 M\acong N&\quad\Longleftrightarrow\quad M^a\cong N^a,\\
 K\ain D^{[u,v]}(R)
 &\quad\Longleftrightarrow\quad K^a\in D^{[u,v]}(R^a)\\
 &\quad\Longleftrightarrow\quad
 \mathfrak aH^q(K)=0\quad(q\notin[u,v]).
 \end{aligned}
\]

We now prove an almost version of
Lemma~\ref{lem:general-completion-amplitude}, with the use of
\cite[Lemma~2.4.2(ii)--(iii)]{GabberRamero}.

\begin{lemma}
\label{lem:general-almost-completion-amplitude}
Let \(R\) be a commutative ring and let
\(\mathfrak a\subset R\) be an ideal satisfying
\(\mathfrak a^2=\mathfrak a\).
Let \((K_s)_{s\geq1}\) be an inverse system in \(D(R)\)
whose transition maps fit into distinguished triangles
\[
 K_1\longrightarrow K_{s+1}\longrightarrow K_s
 \longrightarrow K_1[1]\qquad(s\geq1).
\]
If \(K_1\ain D^{[0,e]}(R)\) for an integer \(e\geq0\), then
\[
 K_s\ain D^{[0,e]}(R)\quad(s\geq1),\qquad
 R\!\varprojlim_s K_s\ain D^{[0,e]}(R).
\]
\end{lemma}

\begin{proof}
By induction after almostification, we get
\(K_s\ain D^{[0,e]}(R)\) for every \(s\geq1\).
By \cite[Lemma~2.4.2(ii)]{GabberRamero}, we have $\varprojlim_sH^q(K_s)=0$ for $q>e$. 
Then the Milnor sequence \eqref{eq:general-lattice-Milnor}, we get
\[
 H^q(R\!\varprojlim_s K_s)\acong R^1\varprojlim_sH^{q-1}(K_s).
\]

Since \(H^{q}(K_1)\acong0\) for $q>e$, the long exact sequences give us
\[
 \mathfrak a\,\operatorname{coker}\!\left(
 H^{q-1}(K_{s+1})\longrightarrow H^{q-1}(K_s)
 \right)=0.
\]
Now we apply \cite[Lemma~2.4.2(iii)]{GabberRamero}, and then we get
\[
 R^1\!\varprojlim_s H^{q-1}(K_s)\acong0.
\]for $q>e$, 
which proves the assertion.
\end{proof}

\subsection{Local computations}
For the rest of this section, let \(\mathcal F\) be a complete unitary
pro-\'etale \(\Qp\)-local system on \(S\).  We choose
\(\mathcal F^+\), \(\mathcal L_s\), and \(\mathcal L_{s,\alpha}\) as in
Definition~\ref{def:complete-unitary-proetale-local-system} and set
\[
 \mathcal F_s:=\mathcal F^+/p^s\mathcal F^+\qquad(s\geq1).
\]
For \(\mathscr B\in\{\widehat{\mathscr O},\mathbb B_J\}\), we write
\[
 \mathcal F^{\mathscr B}
 :=\mathcal F\widehat\otimes_{\Qp}\mathscr B
 \simeq\mathcal F\otimes_{\Qp}^{L,\square}\mathscr B.
\]
The integral models below compute the completed tensor product;
Lemmas~\ref{lem:structural-local-completion}
and~\ref{lem:period-local-completion} identify it with the derived
solid tensor product.

Recall that $S$ is connected. We fix a geometric point $\bar s$ and put $\Gamma=\pi_1^{\et}(S,\bar s)$.
We consider the universal pro-finite-\'etale cover $\widetilde{S}\longrightarrow S$, and denote by
\begin{gather*}
 M_s=(\mathcal L_s)_{\bar s},\qquad
 \Lambda=\varprojlim_s M_s,\qquad V=\Lambda[1/p].
\end{gather*}
By taking stalks of the exact reduction sequences of \(\mathcal F^+\), we get the exact sequence
\[
 0\longrightarrow M_t\xrightarrow{p^s}M_{s+t}
 \longrightarrow M_s\longrightarrow0
 \qquad(s,t\geq1).
\]
Taking \(\varprojlim_t\) along the surjections
\(M_{t+1}\to M_t\), we get
\begin{equation}
 0\longrightarrow\Lambda\xrightarrow{p^s}\Lambda
 \longrightarrow M_s\longrightarrow0,\qquad
 \Lambda\cong\varprojlim_s\Lambda/p^s.
\end{equation} Thus we have the presentation
\(\Lambda\cong(\bigoplus_{b\in B}\Zp)^\wedge_p\) 
by lifting an \(\mathbf F_p\)-basis of \(M_1\).
Therefore \(V\) is Banach with continuous \(\Gamma\)-action and
stable unit ball \(\Lambda\).  On $\widetilde{S}$, all finite local systems are trivial,
and hence we have the following
\begin{equation}\label{eq:general-geometric-trivialization}
 \mathcal F_s|_{\widetilde{S}}\cong\underline{M_s},\qquad
 \mathcal F^+|_{\widetilde{S}}\cong
 \varprojlim_s\underline{M_s}=\underline\Lambda,\qquad
 \mathcal F|_{\widetilde{S}}\cong\underline V.
\end{equation}
Here \(\underline M\) denotes the sheaf of continuous \(M\)-valued functions.
The completed realization of \(\widetilde{S}\) is an affinoid perfectoid
\cite[Theorem~4.9 and its proof]{ScholzeRigid}. Pullback along $\widetilde{S}\to S$ gives us
a basis of affinoid perfectoids on which the sheaves $\cF_s, \cF^+, \cF$ are trivial.

\subsubsection{Tensoring with the structure sheaf}
We first compute the value of $\cF^{\widehat{\mathscr{O}}}(W)$ on an affinoid perfectoid basis of $S_{\proet}$.
\begin{definition}
\label{def:general-structural-coefficients}
We define the following sheaves on $S_{\proet}$ related to $\cF^{\widehat{\mathscr{O}}}$: 
\[
 \mathcal F_s^{\widehat{\mathscr O},+}
 :=\mathcal F_s\otimes_{\Zp/p^s}\widehat{\mathscr O}_S^+/p^s,
 \qquad
 \mathcal F^{\widehat{\mathscr O},+}
 :=\varprojlim_s\mathcal F_s^{\widehat{\mathscr O},+}.
\]
In this subsection, we use the almost mathematics with respect to
\((\cO_C,\mathfrak m_C)\).
\end{definition}

\begin{lemma}
\label{lem:structural-perfectoid-computation}
For every affinoid perfectoid \(W\to S\),
\begin{gather*}
 R\Gamma_{\proet}(W,\mathcal F_s^{\widehat{\mathscr O},+}),\ 
 R\Gamma_{\proet}(W,\mathcal F^{\widehat{\mathscr O},+})
 \ain D^{[0,0]},\\
 R\Gamma_{\proet}(W,\mathcal F^{\widehat{\mathscr O},+})
 \acong R\!\varprojlim_s
 R\Gamma_{\proet}(W,\mathcal F_s^{\widehat{\mathscr O},+}).
\end{gather*}
Furthermore, if \(W=\operatorname{Spa}(R,R^+)\) trivializes
\(\mathcal F^+\) as in \eqref{eq:general-geometric-trivialization},
write
\(\Lambda\cong(\bigoplus_{b\in B}\Zp)^\wedge_p\).  Then
\[
 R\Gamma_{\proet}(W,\mathcal F_s^{\widehat{\mathscr O},+})
 \acong\left(\bigoplus_{b\in B}R^+/p^s\right)[0],\qquad
 R\Gamma_{\proet}(W,\mathcal F^{\widehat{\mathscr O},+})
 \acong\left(\bigoplus_{b\in B}R^+\right)^\wedge_p[0].
\]
Therefore $\cF^{\widehat{\mathscr{O}}}(W)=\mathcal F^{\widehat{\mathscr O},+}(W)[1/p]
 \cong V\widehat\otimes_{\Qp}R
 \cong V\otimes_{\Qp}^{L,\square}R.$
\end{lemma}

\begin{proof}
On the affinoid perfectoid basis trivializing \(\mathcal F^+\), we have
\[
 \mathcal F_s^{\widehat{\mathscr O},+}|_W
 \simeq\underline{\Lambda/p^s}
       \otimes_{\Zp/p^s}\widehat{\mathscr O}_W^+/p^s,
 \qquad \Lambda/p^s\cong\bigoplus_{b\in B}\Zp/p^s.
\]Using almost acyclicity and compatibility with filtered colimits
\cite[Proposition~3.12 and Lemma~4.10]{ScholzeRigid}, we get
\[
 R\Gamma_{\proet}(W,\mathcal F_s^{\widehat{\mathscr O},+})\simeq\bigoplus_{b\in B}R\Gamma_{\proet}(W,\widehat{\mathscr O}_W^+/p^s)
 \acong\left(\bigoplus_{b\in B}R^+/p^s\right)[0].
\]
Thus 
\[
 \mathcal F^{\widehat{\mathscr O},+}(W)[1/p]
 \cong\left(\bigoplus_{b\in B}R^+\right)^\wedge_p[1/p]
 \cong V\widehat\otimes_{\Qp}R
 \cong V\otimes_{\Qp}^{L,\square}R.
\]
The last identification uses Banach solid flatness and
\cite[Corollary~A.28 and Proposition~A.68]{BoscoDrinfeld}.

For an arbitrary affinoid perfectoid \(W\to S\), we can use
\[
 \mathcal F_1|_W
 =\varinjlim_\alpha(\nu^*\mathcal L_{1,\alpha})|_W.
\]
For each finite local system $\cL_{1,\alpha}$, we have almost purity over perfectoids
\[
 R\Gamma_{\proet}
 (W,\nu^*\mathcal L_{1,\alpha}\otimes_{\mathbf F_p}
                       \widehat{\mathscr O}_W^+/p)
 \ain D^{[0,0]}
\]
by \cite[Lemma~4.12]{ScholzeRigid}.  Taking filtered colimits and using d\'evissage, we obtain, for all $s\ge0$
\[R\Gamma_{\proet}(W,\mathcal F_s^{\widehat{\mathscr O},+})\ain D^{[0,0]}.\]
The long exact sequences then show that the transition maps
on sections are almost surjective for every such \(W\).
The almost form of \cite[Lemma~3.18]{ScholzeRigid}, applied on
the affinoid perfectoid basis, and
Lemma~\ref{lem:general-almost-completion-amplitude} with \(e=0\)
therefore give us
\[
 R\Gamma_{\proet}(W,\mathcal F^{\widehat{\mathscr O},+})
 \acong R\!\varprojlim_s
 R\Gamma_{\proet}(W,\mathcal F_s^{\widehat{\mathscr O},+})
 \ain D^{[0,0]}.
\]
\end{proof}

As with the sheaf $\cF^+$, the sheaf $\cF^{\widehat{\mathscr{O}}, +}$ satisfies the following conditions:
\begin{lemma}
\label{lem:structural-local-completion}
On \(S_{\proet}\), these sheaves above satisfy the relations
\[
 0\longrightarrow\mathcal F_1^{\widehat{\mathscr O},+}
 \xrightarrow{p^s}\mathcal F_{s+1}^{\widehat{\mathscr O},+}
 \longrightarrow\mathcal F_s^{\widehat{\mathscr O},+}
 \longrightarrow0,
\]
\[\mathcal F^{\widehat{\mathscr O},+}
 \acong R\!\varprojlim_s\mathcal F_s^{\widehat{\mathscr O},+},\quad
 \mathcal F^{\widehat{\mathscr O},+}
 \otimes_{\Zp}^L\Zp/p^s
 \acong\mathcal F_s^{\widehat{\mathscr O},+},\quad
 \mathcal F^{\widehat{\mathscr O},+}[1/p]
 \simeq\mathcal F^{\widehat{\mathscr O}}.\]
\end{lemma}

\begin{proof}
      Since $\Lambda$ and $\widehat{\mathscr O}^+$ are both $\Zp$-flat, the first sequence is exact.
      The remaining identities of sheaves follow from the local computation in the previous lemma. We write down the details below. 
      Let $W=\mathrm{Spa}(R, R^+)$ be an affinoid perfectoid trivializing $\cF^+$. Then
      \[\cF_s^{\widehat{\mathscr{O}},+}(W)\acong\bigoplus_{b\in B}(R^+/p^s),\quad R^1\varprojlim_s\cF_s^{\widehat{\mathscr{O}},+}(W)\acong0\quad\Rightarrow\quad R\varprojlim_s\cF_s^{\widehat{\mathscr{O}}, +}\acong\varprojlim_{s}\cF_s^{\widehat{\mathscr{O}}, +}=\cF^{\widehat{\mathscr{O}},+}.\]
      Similarly, \[\begin{aligned}
      \cF^{\widehat{\mathscr{O}},+}(W)\otimes_{\Zp}^L\Zp/p^s\acong(\bigoplus_{b\in B}R^+)^{\wedge}_p\otimes_{\Zp}\Zp/p^s\cong\bigoplus_{b\in B}R^+/p^s\acong\cF_s^{\widehat{\mathscr{O}},+}(W),\\ 
      \cF^{\widehat{\mathscr{O}},+}(W)[1/p]\simeq(\bigoplus_{b\in B}R^+)[1/p]\simeq V\otimes_{\Qp}^{L,\square}R\simeq\cF^{\widehat{\mathscr{O}}}(W).
      \end{aligned}\]
      These two identities prove the desired equations\[\mathcal F^{\widehat{\mathscr O},+}
 \otimes_{\Zp}^L\Zp/p^s
 \acong\mathcal F_s^{\widehat{\mathscr O},+},\quad
 \mathcal F^{\widehat{\mathscr O},+}[1/p]
 \simeq\mathcal F^{\widehat{\mathscr O}}.\]
\end{proof}

\subsubsection{Tensoring with the period sheaves}
Similarly, we can compute $\cF^{\mathbb{B}_J}(W)$ on an affinoid perfectoid basis of $S_{\proet}$.
\begin{definition}
\label{def:Fontaine-Ainf-sheaf}
We define the following sheaves on \(S_{\proet}\):
\[
 \widehat{\mathscr O}^{\flat,+}
 :=\varprojlim_{x\mapsto x^p}\widehat{\mathscr O}^{+}/p,
 \qquad
 \mathbb A_{\inf}
 :=\mathrm W\!\left(\widehat{\mathscr O}^{\flat,+}\right),
\]
where \(\mathrm W(-)\) denotes the functor of \(p\)-typical Witt
vectors; see \cite[Definition~6.1]{ScholzeRigid}.
For an affinoid perfectoid test object
\(W=\operatorname{Spa}(R,R^+)\to S\), they have values
\[
 \widehat{\mathscr O}^{\flat,+}(W)=R^{\flat,+},
 \qquad
 \mathbb A_{\inf}(W)=\mathrm W(R^{\flat,+}).
\]
We write \([x]\) for the Teichm\"uller lift of a section
\(x\) of \(\widehat{\mathscr O}^{\flat,+}\) and denote 
\(A_{\inf}(\cO_{C^\flat}):=\mathrm W(\cO_{C^\flat})\).
\end{definition}

Fix a compact interval \(J=[r_0,r_1]\subset(0,\infty)\)
with rational endpoints.  Choose a pseudo-uniformizer
\(\rho\in C^\flat\) with compatible roots, and put
\begin{equation*}
 \bar{a}_J=\rho^{1/r_1},\qquad
 \bar{b}_J=\rho^{1/r_0},\qquad a_J=[\bar{a}_J].
\end{equation*}
Following \cite[Definition~4.4]{BoscoDrinfeld}, we define the sheaf of rings
\begin{equation*}
 \mathbb A_J
 :=\left(\mathbb A_{\inf}
       \left[\frac p{[\bar{a}_J]},\frac{[\bar{b}_J]}p\right]
       \right)^\wedge_p,
 \qquad \mathbb B_J=\mathbb A_J[1/p].
\end{equation*}
The $\mathbb{A}_J$ is the integral model \(\mathbb B_J^+\) of
\cite[Theorem~2.2.4 and its proof]{DRC}.
For an integer \(N\geq r_1/r_0\), we have relations
\begin{equation}\label{eq:period-cofinal-filtrations}
 a_J^N\mathbb A_J\subset p\mathbb A_J\subset a_J\mathbb A_J,
 \qquad
 \mathbb A_J[1/a_J]=\mathbb A_J[1/p].
\end{equation}
Indeed, \(p=a_J(p/a_J)\) and
\(a_J^N=p([\bar{b}_J]/p)[\bar{a}_J^N/\bar{b}_J]\).
Thus the \(a_J\)-adic and \(p\)-adic filtrations are cofinal.
By \cite[Theorem~2.2.4 and Corollary~2.2.11]{DRC},
\(\mathbb A_J\) is \(a_J\)-torsion-free and
\[
 \mathbb A_J/a_J
 \cong
 (\widehat{\mathscr O}^{\flat,+}/\bar{a}_J)
 [T_0,T_1]/(T_0T_1-\bar{b}_J/\bar{a}_J)
\cong\bigoplus_{k\geq1}
       \widehat{\mathscr O}^{\flat,+}/\bar{a}_J
 \quad\text{as }\mathbb A_{\inf}\text{-modules}.
 \label{eq:period-annulus-free-reduction}
\]
In this subsection, we therefore use almost mathematics with respect to
\begin{equation*}
 \left(A_{\inf}(\cO_{C^\flat}),\mathfrak m_{\inf}\right),
 \qquad
 \mathfrak m_{\inf}
   =\bigl([\rho^{1/p^k}]:k\geq0\bigr).
\end{equation*}
The ideal becomes the unit ideal after inverting $p$, thus in $\mathbb{B}_J$, we have
\(
 \mathfrak m_{\inf}\mathbb B_J(C)=\mathbb B_J(C).
\)
Thus Lemma~\ref{lem:general-almost-completion-amplitude} applies
over \(\mathbb A_J(C)\) with the extended almost ideal $m_{\inf}\mathbb A_J(C)$.

\begin{definition}
\label{def:general-period-coefficients}
On \(S_{\proet}\), set
\[
 \mathcal F_s^{\mathbb B_J,+}
 :=\mathcal F^+\otimes_{\Zp}(\mathbb A_J/a_J^s),
 \qquad
 \mathcal F^{\mathbb B_J,+}
 :=\varprojlim_s\mathcal F_s^{\mathbb B_J,+}.
\]
The ordinary sheaf tensor product depends only on
\(\mathcal F_s\), since \(p^s\mathbb A_J\subset a_J^s\mathbb A_J\).
\end{definition}

\begin{lemma}
\label{lem:period-integral-perfectoid}
For every affinoid perfectoid \(W\to S\) and every \(s\geq1\),
\begin{equation}\label{eq:period-perfectoid-torsion-vanishing}
 \begin{aligned}
 R\Gamma_{\proet}(W,\mathcal F_s^{\mathbb B_J,+})
 &\ain D^{[0,0]}\qquad(s\geq1),\\
 R\Gamma_{\proet}(W,\mathcal F^{\mathbb B_J,+})
 &\acong
 \left(\varprojlim_s
 H^0(W,\mathcal F_s^{\mathbb B_J,+})\right)[0].
 \end{aligned}
\end{equation}
\end{lemma}

\begin{proof}
Let $W=\operatorname{Spa}(R,R^+)\to S$ be an affinoid perfectoid.
Since $p=0$ on $\mathbb A_J/a_J$,
\[
 \mathcal F_1^{\mathbb B_J,+}
 =\mathcal F_1\otimes_{\mathbf F_p}
                 (\mathbb A_J/a_J),\qquad
 \mathcal F_1
 =\varinjlim_\alpha\nu^*\mathcal L_{1,\alpha}.
\]
Fix \(\alpha\) and put
\(\mathcal L=(\nu^*\mathcal L_{1,\alpha})|_W\).
Choose a finite
\'etale \(\Delta\)-torsor \(S_\alpha\to S\) trivializing
\(\mathcal L_{1,\alpha}\).  Pulling back to $W$, we denote
\[
 W_\alpha=W\times_SS_\alpha=\operatorname{Spa}(R_\alpha,R^+_\alpha)
 \longrightarrow W,\qquad
 \mathcal L|_{W_\alpha}\cong\underline N.
\]
Then \(W_\alpha\) is an affinoid perfectoid.  Put
\(
 R_0=R^{\flat,+}/\bar{a}_J,
 S_0=R^{\flat,+}_\alpha/\bar{a}_J.\) By \cite[Theorems~7.9 and~7.12]{ScholzePerfectoid}
, we have tilting equivalence and almost purity 
\[
 W_{\et}\simeq W^\flat_{\et},\qquad
 S_0/R_0\text{ almost finite \'etale and almost faithfully flat},
 \qquad
 S_0^a\otimes_{R_0^a}S_0^a
 \xrightarrow{\sim}\prod_{\gamma\in\Delta}S_0^a,
\]
where $x\otimes y\mapsto(x\gamma(y))_\gamma$.
In characteristic $p$, we use the almost mathematics with respect to
$(\cO_{C^\flat},\mathfrak m_{C^\flat})$.
We have almost acyclicity theorem by \cite[Proposition~7.13]{ScholzePerfectoid} and
\cite[Proposition~3.12 and Lemma~4.10]{ScholzeRigid} that
\[
 R\Gamma_{\proet}
 (W_\alpha,\widehat{\mathscr O}^{\flat,+}/\bar{a}_J)
 \acong S_0[0].
\]
The right $\Delta$-torsor $W_\alpha\to W$ has \v Cech nerve
\[
 W_\alpha\times\Delta^q\xrightarrow{\sim}
 W_\alpha^{\times q}:=\underbrace{W_\alpha\times_W\cdots\times_W W_\alpha}_{q+1\text{ factors}}
\]
The $\Delta$-action on $N\otimes_{\mathbf F_p}S_0$ is given by
$g(n\otimes s)=gn\otimes r_g^*s$, where $r_g(w)=wg$.
Using the trivialization along $W_\alpha^{\times q}\to W_\alpha$, 
we identify the sheaf cochain complex with the group cochain complex:
\[
 \check C^q:=
 H^0\!\left(W_\alpha^{\times q},\mathcal L\otimes_{\mathbf F_p}
       \widehat{\mathscr O}^{\flat,+}/\bar{a}_J\right)
 \acong\operatorname{Map}(\Delta^q,N\otimes_{\mathbf F_p}S_0).
\]
Using almost acyclicity on every perfectoid $W_\alpha^{\times q}$ again, we obtain
\[
 R\Gamma_{\proet}\!\left(W,
   \mathcal L\otimes_{\mathbf F_p}
       \widehat{\mathscr O}^{\flat,+}/\bar{a}_J\right)
 \acong\check C^\bullet
 \acong R\Gamma(\Delta,N\otimes_{\mathbf F_p}S_0).
\]
By \cite[proof of Lemma~4.12]{ScholzeRigid}, 
the almost faithfully flat descent for $S_0/R_0$
then gives us
\[
 R\Gamma(\Delta,N\otimes_{\mathbf F_p}S_0)\acong N\otimes_{\mathbf F_p}S_0^\Delta[0],\qquad
\]
By \cite[Proposition~3.12 and Lemma~4.10]{ScholzeRigid}, 
we can take filtered colimit for the pro-\'etale cohomology on $W$.  
The decomposition of $\mathbb{A}_J/a_J$ allows us to compute 
\begin{align*}
 R\Gamma_{\proet}(W,\mathcal F_1^{\mathbb B_J,+})
 &\simeq\varinjlim_\alpha\bigoplus_{\ell\ge1}
 R\Gamma_{\proet}\!\left(W,
 \nu^*\mathcal L_{1,\alpha}\otimes_{\mathbf F_p}
       \widehat{\mathscr O}^{\flat,+}/\bar{a}_J\right)\\
 &\acong\left(\varinjlim_\alpha\bigoplus_{\ell\ge1}
 N\otimes_{\mathbf F_p}S_0^\Delta\right)[0].
\end{align*}
In particular, $H^q_{\proet}(W,\mathcal F_1^{\mathbb B_J,+})\acong0$ for $q>0$.
Since $\mathcal F^+$ is flat over $\Zp$ and $\mathbb{A}_J$ is
$a_J$-torsion-free,
\[
 0\longrightarrow\mathcal F_1^{\mathbb B_J,+}
 \xrightarrow{a_J^s}\mathcal F_{s+1}^{\mathbb B_J,+}
 \longrightarrow\mathcal F_s^{\mathbb B_J,+}\longrightarrow0.
\]
Apply Lemma~\ref{lem:general-almost-completion-amplitude} to the case
\(e=0\) and \(K_s=R\Gamma_{\proet}(W,\mathcal F_s^{\mathbb B_J,+})\), and then we can get
\[
 R^1\!\varprojlim_s H^0(W,\mathcal F_s^{\mathbb B_J,+})\acong0,\qquad
 R\!\varprojlim_sR\Gamma_{\proet}(W,\mathcal F_s^{\mathbb B_J,+})
 \acong\left(\varprojlim_sH^0(W,\mathcal F_s^{\mathbb B_J,+})\right)[0].
\]
Finally we apply the almost form of \cite[Lemma~3.18]{ScholzeRigid}
on the affinoid perfectoid basis. Consequently 
\[
 \mathcal F^{\mathbb B_J,+}
 \acong R\!\varprojlim_s\mathcal F_s^{\mathbb B_J,+},\qquad
 R\Gamma_{\proet}(W,\mathcal F^{\mathbb B_J,+})
 \acong R\!\varprojlim_sR\Gamma_{\proet}(W,\mathcal F_s^{\mathbb B_J,+}).
\]
\end{proof}

\begin{lemma}
\label{lem:period-local-completion}
On \(S_{\proet}\), the completed period sheaf satisfies
\[
 \mathcal F^{\mathbb B_J,+}
 \acong R\!\varprojlim_s\mathcal F_s^{\mathbb B_J,+},
 \quad
 \mathcal F^{\mathbb B_J,+}[1/p]
 \simeq\mathcal F^{\mathbb B_J}.
\]
\end{lemma}

\begin{proof}
The first equation of derived completeness was proved in Lemma~\ref{lem:period-integral-perfectoid}.
To identify the second equation of localization, we use the basis $W\to S$ with
\(\mathcal F^+|_W\cong\underline\Lambda\) from
\eqref{eq:general-geometric-trivialization}.
By \cite[Proposition~4.7]{BoscoDrinfeld}
and $a_J$-torsion-freeness,
\[
 R\Gamma_{\proet}(W,\mathbb A_J)\acong \mathbb A_J(W)[0],
 \qquad
 R\Gamma_{\proet}(W,\mathbb A_J/a_J^s)
 \acong(\mathbb A_J(W)/a_J^s)[0].
\]
Choose a Banach basis
$\Lambda\cong(\bigoplus_{b\in B}\Zp)^\wedge_p$.
Then
\[
 \mathcal F_s^{\mathbb B_J,+}|_W
 \cong\bigoplus_{b\in B}(\mathbb A_J/a_J^s)|_W,\qquad
 \mathcal F^{\mathbb B_J,+}(W)
 \acong\varprojlim_s\bigoplus_{b\in B}\mathbb A_J(W)/a_J^s.
\]
Choose an integer $c\ge r_1/r_0$.  For $m\ge1$,
we have inclusions
$a_J^{cm}\mathbb A_J(W)\subset p^m\mathbb A_J(W)\subset a_J^m\mathbb A_J(W)$.
Since the almost ideal becomes the unit ideal after inverting
$p$, we obtain
\begin{align*}
 \mathcal F^{\mathbb B_J,+}(W)[1/p]
 &\cong\left(\bigoplus_{b\in B}\mathbb A_J(W)\right)^\wedge_{a_J}[1/p]\cong\left(\bigoplus_{b\in B}\mathbb A_J(W)\right)^\wedge_p[1/p]\\
 &\cong V\widehat\otimes_{\Qp}\mathbb B_J(W)
  \cong V\otimes_{\Qp}^{L,\square}\mathbb B_J(W).
\end{align*}
Banach spaces are solid flat and thus by
\cite[Corollary~A.28 and Proposition~A.68]{BoscoDrinfeld}
we get the last isomorphism.  The natural tensor maps descend to
$\mathcal F^{\mathbb B_J,+}[1/p]\simeq\mathcal F^{\mathbb B_J}$.
\end{proof}

\begin{lemma}
\label{lem:general-geometric-perfectoid-parameters}
\label{lem:period-coefficient-acyclicity}
For every affinoid perfectoid \(W\to S\) and
\(\mathscr B\in\{\widehat{\mathscr O},\mathbb B_J\}\), the space
\(H^0(W,\mathcal F^{\mathscr B})=\mathcal F^{\mathscr B}(W)\)
is a Banach space, and for every profinite set \(Z\) there is a natural
equivalence
\begin{equation}\label{eq:general-perfectoid-parameters}
 R\Gamma_{\proet}(W\times Z,\mathcal F^{\mathscr B})
 \simeq C(Z,\mathcal F^{\mathscr B}(W))[0].
\end{equation}
The identifications commute with pullback and pro-\'etale descent.
\end{lemma}

\begin{proof}
Since \(W\times Z\) is affinoid perfectoid, combining
Lemmas~\ref{lem:structural-perfectoid-computation},
\ref{lem:structural-local-completion},
\ref{lem:period-integral-perfectoid}, and
\ref{lem:period-local-completion}, followed by
Lemma~\ref{lem:general-almost-completion-amplitude} with \(e=0\), we obtain the concentration
\[
 R\Gamma_{\proet}(W\times Z,\mathcal F^{\mathscr B})
 \in D^{[0,0]}.
\]
Use the data
\((\widetilde{S},\Gamma,\Lambda,V)\) as in
\eqref{eq:general-geometric-trivialization}.
Set \(\widetilde{W}=W\times_S\widetilde{S}\), an affinoid perfectoid
pro-finite-\'etale \(\Gamma\)-torsor over $W$. The cohomology descent is concentrated in degree zero 
and thus gives
\[
 \mathcal F^{\mathscr B}(W)
 \cong\left(V\widehat\otimes_{\Qp}\mathscr B(\widetilde{W})\right)^\Gamma.
\]
These invariants form a closed Banach subspace.
For affinoid perfectoid \(W\),
\[
 \widehat{\mathscr O}(W\times Z)
 =C(Z,\widehat{\mathscr O}(W)),\qquad
 \mathbb B_J(W\times Z)=C(Z,\mathbb B_J(W)).
\]
We obtain the first identity by completing the rings for
finite disjoint unions of \(Z\). The second identity follows from
\cite[Corollary~4.10]{BoscoDrinfeld}.
For Banach \(\Qp\)-spaces \(V,A\),
\cite[Proposition~A.55(i), Corollary~A.61, and
Proposition~A.68]{BoscoDrinfeld} give us the equation
\[
 V\widehat\otimes_{\Qp}C(Z,A)
 \cong C(Z,V\widehat\otimes_{\Qp}A).
\]
Since \(\Gamma\) acts trivially on \(Z\),
\begin{align*}
 H^0(W\times Z,\mathcal F^{\mathscr B})
 &\cong
 \left(V\widehat\otimes C(Z,\mathscr B(\widetilde{W}))\right)^\Gamma\\
 &\cong
 C\!\left(Z,
   \left(V\widehat\otimes\mathscr B(\widetilde{W})\right)^\Gamma\right)
 =C(Z,\mathcal F^{\mathscr B}(W)).
\end{align*}
\end{proof}

\subsection{The general affinoid bound}
\label{subsec:general-affinoid-bound}

\begin{lemma}
\label{lem:general-Banach-group-dimension}
Let \(\Gamma\) be a profinite group with
\(\operatorname{cd}_p(\Gamma)\leq d\).
For every continuous Banach \(\Qp\)-representation \(M\) of $\Gamma$,
we have 
\[
 R\Gamma_{\mathrm{cont}}(\Gamma,M)\in D^{[0,d]}.
\]
The assertion holds in the solid category as well.
\end{lemma}
\begin{proof}
Choose a complete \(\Gamma\)-stable open bounded lattice
\(M^+\subset M\), using the invariant norm
\(\|m\|_\Gamma=\sup_{\gamma\in\Gamma}\|\gamma m\|\).
Put \(M_s=M^+/p^s\), with its discrete topology.
For every \(q\geq0\), since $\Gamma$ is a compact group, we have the following relations
\begin{gather*}
 C_{\mathrm{cont}}(\Gamma^q,M^+)
 =\varprojlim_s C_{\mathrm{cont}}(\Gamma^q,M_s),\quad
 C_{\mathrm{cont}}(\Gamma^q,M)
 =C_{\mathrm{cont}}(\Gamma^q,M^+)[1/p],\\
 C_{\mathrm{cont}}(\Gamma^q,M_{s+1})
 \twoheadrightarrow C_{\mathrm{cont}}(\Gamma^q,M_s).
\end{gather*}
Surjectivity follows by lifting finitely many values.  Thus
\[
 R\Gamma_{\mathrm{cont}}(\Gamma,M^+)
 \simeq R\!\varprojlim_s R\Gamma_{\mathrm{cont}}(\Gamma,M_s).
\]
The exact sequences
\[
 0\longrightarrow M_1\xrightarrow{p^s}M_{s+1}
 \longrightarrow M_s\longrightarrow0
\]
give us, in the top degree,
\[
 H^d_{\mathrm{cont}}(\Gamma,M_{s+1})
 \longrightarrow H^d_{\mathrm{cont}}(\Gamma,M_s)
 \longrightarrow H^{d+1}_{\mathrm{cont}}(\Gamma,M_1)=0,
\]
since \(M_1\) is discrete, killed by \(p\), and
\(\operatorname{cd}_p\Gamma\leq d\).
Thus the top-degree transitions are surjective; the integral
and rational bounds follow from
Lemma~\ref{lem:general-completion-amplitude} and inverting $p$.
For an extremally disconnected profinite set \(Z\),
the exponential law identifies evaluated cochains:
\[
 C(Z,C(\Gamma^q,M))\cong C(\Gamma^q,C(Z,M)).
\]
Apply the same proof to the Banach representation \(C(Z,M)\).
These exact evaluations prove the corresponding bound in the solid category.
\end{proof}

We now identify the pro-\'etale cohomology of $\cF^{\mathscr B}$
with the continuous cohomology of the profinite group \(\Gamma\). 
Recall that $\widetilde{S}\to S$ is the universal pro-finite-\'etale cover with Galois group $\Gamma$, and $\cF^+|_{\widetilde{S}}\cong\underline{\Lambda}$. 
As in the preceding subsection, we use the following almost settings:
\begin{itemize}
\item If \(\mathscr{B}=\widehat{\mathscr{O}}\), the almost mathematics is with respect to the pair
\((\cO_C, \mathfrak{m}_C)\); put \(\mathscr B^+=\widehat{\mathscr O}^+\) and \(t=p\).
\item If \(\mathscr{B}=\mathbb{B}_J\), the almost mathematics is with respect to the pair \((A_{\inf}(\cO_{C^\flat}), \mathfrak{m}_{\inf})\); put \(\mathscr B^+=\mathbb A_J\) and \(t=a_J\).   
\end{itemize} 
\begin{proposition}
\label{prop:geometric-almost-Cartan-Leray}
\label{prop:structural-almost-Cartan-Leray}
Let $\cF$ be a complete unitarity pro-\'etale local system with lattice $\cF^+$. 
Let $\cF^{\mathscr{B}, +}$, $\cF^{\mathscr{B}, +}_s$ be as in Definitions~\ref{def:general-structural-coefficients} or \ref{def:general-period-coefficients}. 
With the diagonal \(\Gamma\)-actions, the \(t\)-adic topology
on the completed objects, and the discrete topology on
the reductions, we have
\begin{align}
 R\Gamma_{\proet}(S,\mathcal F_s^{\mathscr B,+})
 &\acong R\Gamma_{\mathrm{cont}}\bigl(\Gamma,
 \Lambda\otimes_{\Zp}(\mathscr B^+(\widetilde{S})/t^s)\bigr)
       \in D^{[0,d]},\\
 R\Gamma_{\proet}(S,\mathcal F^{\mathscr B,+})
 &\acong R\Gamma_{\mathrm{cont}}\bigl(\Gamma,
 \Lambda\otimes_{\Zp}^{\square}\mathscr B^+(\widetilde{S})\bigr)
       \in D^{[0,d]}.
\end{align}
\end{proposition}
\begin{proof}
The ring \(\mathscr B^+(\widetilde{S})\) is \(t\)-torsion-free and
\(t\)-adically complete.  For \(\mathbb A_J\), this follows from
its integral construction and the relation
\eqref{eq:period-cofinal-filtrations}.
By choosing a basis of $\Lambda/p$ over $\mathbb{F}_p$, we can express $\Lambda$ as \(\Lambda\cong(\bigoplus_{b\in B}\Zp)^\wedge_p\).
Since \(p^s\mathscr B^+(\widetilde{S})\subset
 t^s\mathscr B^+(\widetilde{S})\), we have
\[
 \Lambda\otimes_{\Zp}(\mathscr B^+(\widetilde{S})/t^s)
 \cong\bigoplus_{b\in B}\mathscr B^+(\widetilde{S})/t^s,\qquad
 \Lambda\otimes_{\Zp}^{\square}\mathscr B^+(\widetilde{S})
 \cong\left(\bigoplus_{b\in B}\mathscr B^+(\widetilde{S})\right)^\wedge_t.
\]

For every profinite set \(Z\), by
\cite[Proposition~4.7(i)--(ii) and Corollary~4.10]{BoscoDrinfeld} we have
\begin{align*}
 \mathscr B^+(\widetilde{S}\times Z)
 &=C(Z,\mathscr B^+(\widetilde{S})),\\
 C(Z,\mathscr B^+(\widetilde{S}))/t^s
 &\cong C(Z,\mathscr B^+(\widetilde{S})/t^s).
\end{align*}
The second identity holds because $Z$ is a compact space, and division by \(t^s\)
on \(t^s\mathscr B^+(\widetilde{S})\) is continuous for the
\(t\)-adic topology.
For \(Z=\Gamma\), pulling back along the action morphism $a\colon\widetilde{S}\times\Gamma\to\widetilde{S}, (s,\gamma)\mapsto \gamma s$,
gives a continuous map 
\[
 \mathscr B^+(\widetilde{S})\xrightarrow{a^*}
 C(\Gamma,\mathscr B^+(\widetilde{S})).
\]
Hence $\Gamma$ acts continuously on \(\mathscr B^+(\widetilde{S})\) 
and its discrete quotients.  The action on \(\Lambda/p^s\)
is also continuous, hence so is that on the tensor product
\(\Lambda\otimes_{\Zp}(\mathscr B^+(\widetilde{S})/t^s)\).

Thus by Lemmas~\ref{lem:structural-perfectoid-computation} or
\ref{lem:period-integral-perfectoid} for the two choice of $\mathscr{B}$, and
\cite[Proposition~4.7]{BoscoDrinfeld}, together with compactness of $Z$, we have an almost isomorphism for each $Z$
\begin{equation}\label{eq:geometric-integral-profinite-parameters}
 R\Gamma_{\proet}(\widetilde{S}\times Z,\mathcal F_s^{\mathscr B,+})
 \acong C\bigl(Z,\Lambda\otimes_{\Zp}
 (\mathscr B^+(\widetilde{S})/t^s)\bigr)[0].
\end{equation}
Apply this to
\(\widetilde{S}^{\times_{S}(r+1)}\cong \widetilde{S}\times\Gamma^r\) and the \v Cech nerve $[\dots\to\widetilde{S}^{\times_S(r+1)}\to\widetilde{S}^{\times_S(r)}\to\dots\to\widetilde{S}\to S]$.
The \v Cech complex can be identified with the continuous cochain of group cohomology:
\[
 R\Gamma_{\proet}(S,\mathcal F_s^{\mathscr B,+})
 \acong C\bigl(\Gamma^\bullet,\Lambda\otimes_{\Zp}
 (\mathscr B^+(\widetilde{S})/t^s)\bigr)= R\Gamma_{\mathrm{cont}}\bigl(\Gamma,
 \Lambda\otimes_{\Zp}(\mathscr B^+(\widetilde{S})/t^s)\bigr)
\]
Since $\Lambda\otimes_{\mathbb{Z}_p}\mathscr{B}^+(\widetilde{S})/t^s=\Lambda/p^s\otimes_{\mathbb{Z}_p/p^s}\mathscr{B}^+(\widetilde{S})/t^s$ 
carries the discrete topology, by Proposition~\ref{prop:affinoid-fundamental-group-dimension} we get
\[
 H^q_{\mathrm{cont}}\bigl(\Gamma,
 \Lambda\otimes_{\Zp}(\mathscr B^+(\widetilde{S})/t^s)\bigr)=0
 \qquad(q>d).
\]

Taking inverse limit to the complex above, we obtain
\begin{align*}
 R\Gamma_{\proet}(S,\mathcal F^{\mathscr B,+})
 &\acong R\!\varprojlim_s R\Gamma_{\mathrm{cont}}\bigl(\Gamma,
 \Lambda\otimes_{\Zp}(\mathscr B^+(\widetilde{S})/t^s)\bigr)\\
 &\simeq R\Gamma_{\mathrm{cont}}\bigl(\Gamma,
 \Lambda\otimes_{\Zp}^{\square}\mathscr B^+(\widetilde{S})\bigr).
\end{align*}
The last equivalence follows by lifting finite images in each
cochain degree.  The \(\Gamma\)-equivariant exact sequences
\[
 0\longrightarrow
 \Lambda\otimes_{\Zp}(\mathscr B^+(\widetilde{S})/t)\xrightarrow{t^s}
 \Lambda\otimes_{\Zp}(\mathscr B^+(\widetilde{S})/t^{s+1})\longrightarrow
 \Lambda\otimes_{\Zp}(\mathscr B^+(\widetilde{S})/t^s)
 \longrightarrow0
\]
then give surjective transition maps on \(H^d\) terms, hence
\[
 R^1\!\varprojlim_s H^d_{\mathrm{cont}}\bigl(\Gamma,
 \Lambda\otimes_{\Zp}(\mathscr B^+(\widetilde{S})/t^s)\bigr)\acong0.
\] The completed bound then follows from 
Lemma~\ref{lem:general-completion-amplitude}.  Lemmas~\ref{lem:structural-local-completion} and
\ref{lem:period-local-completion} identify the complexes after inverting $p$.
\end{proof}

We can now prove the main theorem with coefficients.

\begin{theorem}
\label{thm:general-unitary-coefficient-amplitude}
Let \(\mathcal F\) be a complete unitary pro-\'etale \(\Qp\)-local system
on \(S\), with a lattice \(\mathcal F^+\).
For every compact rational interval \(J\subset(0,\infty)\)
and every \(\mathscr B\in\{\Qp,\widehat{\mathscr O},\mathbb B_J\}\),
we have
\begin{equation}\label{eq:general-three-coefficient-amplitude}
 R\Gamma_{\proet}(S,\mathcal F^+),\ 
 R\Gamma_{\proet}(S,\mathcal F\widehat{\otimes}_{\Qp}\mathscr B)\in D^{[0,d]}.
\end{equation}
The completed tensor products agree with
\(\mathcal F\otimes_{\Qp}^{L,\square}\mathscr B\), and the bounds hold in the solid category.
\end{theorem}
\begin{proof}
The case $\mathscr{B}=\mathbb{Q}_p$ is proved in Theorem~\ref{thm:general-local-lattice-amplitude}.
For \(\mathscr B\in\{\widehat{\mathscr O},\mathbb B_J\}\), we apply Proposition~\ref{prop:geometric-almost-Cartan-Leray}.
The almost isomorphism in Proposition~\ref{prop:geometric-almost-Cartan-Leray} and the compactness of $\Gamma$ allow us to identify
\begin{equation}\label{eq:general-geometric-Cartan-Leray}
 \begin{aligned}
 R\Gamma_{\proet}(S,\mathcal F^{\mathscr B})
 &\simeq R\Gamma_{\mathrm{cont}}\bigl(\Gamma,
 \Lambda\otimes_{\Zp}^{\square}\mathscr B^+(\widetilde{S})\bigr)[1/p]\\
 &\simeq R\Gamma_{\mathrm{cont}}\bigl(\Gamma,
 V\widehat\otimes_{\Qp}\mathscr B(\widetilde{S})\bigr)
 \in D^{[0,d]}.
 \end{aligned}
\end{equation}
For an extremally disconnected profinite set \(Z\),
Lemma~\ref{lem:general-geometric-perfectoid-parameters} identifies
the evaluated descent complex with
\[
 R\Gamma_{\mathrm{cont}}\bigl(\Gamma,
 C(Z,V\widehat\otimes_{\Qp}\mathscr B(\widetilde{S}))\bigr).
\]
Then by Lemma~\ref{lem:general-Banach-group-dimension} we deduce the bound in the solid category.
\end{proof}


\section{Cohomological amplitude of Jacquet--Langlands}
\label{sec:cohomological-amplitude}

Let \(F/\Qp\) and \(E/\Qp\) be finite extensions, with
\(E\) sufficiently large, and let \(n\geq2\).  Put
\begin{equation*}
 C=\widehat{\overline F},\qquad
 G=\mathrm{GL}_n(F),\qquad
 X=\Period.
\end{equation*}
Let \(D/F\) be the central division algebra of invariant \(1/n\).
The infinite-level Lubin--Tate space \(\cM\) has a
\(D^\times\)-equivariant pro-\'etale \(G\)-torsor
\begin{equation*}
 \pi=\pi_{\mathrm{GH}}:\cM\longrightarrow X,\qquad
 \cM/G\simeq X,
\end{equation*}
with commuting \(G\)- and \(D^\times\)-actions
\cite[Proposition~4.4.4]{DRC}.  For a closed subgroup
\(H\subset G\), we write
\begin{equation*}
 \cM_H=\cM/H,\qquad \pi_H:\cM_H\longrightarrow X.
\end{equation*}

For a right \(H\)-space \(T\) and a topological left
\(H\)-module \(\Pi\), set
\begin{equation*}
 C_H(|T|,\Pi)=\{f:|T|\to \Pi\text{ continuous}:\
                 f(th)=h^{-1}f(t)\}.
\end{equation*}
Following \cite{DRC}, for a topological \(G\)-module \(\Pi\), set
\begin{align*}
 \cF_\Pi(W)&=C_G\bigl(|W\times_X\cM|,\Pi\bigr),\qquad W\in X_{\proet},\\
 \JL(\Pi)&=R\Gamma_{\proet}(X,\cF_\Pi),\qquad
 \JL^i(\Pi)=H^i\JL(\Pi).
\end{align*}
The \(D^\times\)-action is \((a\cdot f)(x)=f(a^{-1}x)\).
We use derived categories of solid representations over the
respective coefficient rings and write \(\otimes_E^{L,\square}\)
for the derived solid tensor product.

For \(T/C\), write \(\nu:T_{\proet}\to T_{\et}\) for the
change of sites and \(\underline N_T\) for the sheaf of
continuous functions valued in a topological module \(N\).

Locally analytic means locally \(\Qp\)-analytic.
The compact-type bounds also apply to locally \(F\)-analytic
actions, since they use only continuity.

For a compact interval \(J\subset(0,\infty)\) with rational
endpoints, let
\(\mathscr B\in\{\Qp,\widehat{\mathscr O},\mathbb B_J\}\).
For a representation \(\Pi\), set
\begin{align*}
 \cF_\Pi^{\mathscr B}
 :=\cF_\Pi\otimes_{\Qp}^{L,\square}\mathscr B,\quad
 &\JL_{\mathscr B}(\Pi)
 :=R\Gamma_{\proet}(X,\cF_\Pi^{\mathscr B}),\\
 \cF_\Pi^{\Qp}=\cF_\Pi,\quad&\JL_{\Qp}(\Pi)=\JL(\Pi).
\end{align*}
We tensor over \(\Qp\) to keep the \(E\)-action without
choosing an embedding \(E\hookrightarrow\mathscr B\).
These sheaves are ordinary \(D^\times\)-equivariant sheaves due to
the solid flatness for Banach \(\Pi\) \cite[\S2.3.1]{DRC}
and from filtered colimits for the LB representations of
Subsection~\ref{subsec:LB-lalg-amplitude}.

\subsection{Rational charts and complete unitary local systems}
\label{subsec:local-unitarity}

Fix \(K=\mathrm{GL}_n(\cO_F)\) and a uniformizer
\(\pi_F\in\cO_F\). We first construct a rational covering $\{U_i\}$ of $X=\mathbb{P}^{n-1}$
by rational affinoids, and each $U_i$ admits a section $\sigma\colon U_i\to M_K$. 

\begin{proposition}
\label{prop:GH-rational-cover}
There are homogeneous coordinates \([y_0:\cdots:y_{n-1}]\)
on \(X\) with
\begin{equation*}
 X=\bigcup_{i=0}^{n-1}U_i,\qquad
 U_i=\{ |y_j|\leq |y_i|\text{ for all }j\},\qquad
 \sigma_i:U_i\longrightarrow\cM_K,\quad
 \pi_K\sigma_i=\operatorname{id}_{U_i}.
\end{equation*}
For \(\varnothing\ne I\subset\{0,\ldots,n-1\}\), put
\(U_I:=\bigcap_{i\in I}U_i\).  Each \(U_I\) is a rational
affinoid domain in every \(U_i\), \(i\in I\), and has the
section \(\sigma_i|_{U_I}\).
\end{proposition}

\begin{proof}
In Gross--Hopkins coordinates \([x_0:\cdots:x_{n-1}]\),
there is a closed polydisc \(\Delta\) in the height-zero
component of \(\cM_K\) for which
\begin{equation*}
 \pi_K|_\Delta:\Delta\xrightarrow{\sim}\Omega_0,
 \qquad
 \Omega_0=\left\{\left|\frac{x_j}{x_0}\right|
           \leq|\pi_F|^{(n-j)/n},\quad 1\leq j<n\right\}.
\end{equation*}
See \cite[Lemma~23.14 and Corollary~23.15]{GrossHopkinsBundles}
and \cite[\S3.1, equations~(3.1.8)--(3.1.10)]{ShethLocallyAnalytic}.
Choose a uniformizer \(\delta\in D^\times\), with
\(\delta^n=\pi_F\), and the cyclic basis in which its
projective action is
\begin{equation*}
 \delta[x_0:x_1:\cdots:x_{n-1}]
   =[x_{n-1}:\pi_F x_0:x_1:\cdots:x_{n-2}].
\end{equation*}
Choose \(t\in C\) with \(t^n=\pi_F\), and set
\begin{equation*}
 y_0=x_0,\qquad y_j=t^{j-n}x_j\quad(1\leq j<n).
\end{equation*}
Then \(\Omega_0=U_0\).  In column-vector notation,
\[
 \mathbf y=S\mathbf x,\qquad
 S=\operatorname{diag}(1,t^{1-n},\ldots,t^{-1}),\qquad
 S^{-1}=\operatorname{diag}(1,t^{n-1},\ldots,t).
\]
The action matrices satisfy
\[
 A_\delta=
 \begin{pmatrix}
  0&0&\cdots&0&1\\
  \pi_F&0&\cdots&0&0\\
  0&1&\ddots&\vdots&\vdots\\
  \vdots&\ddots&\ddots&0&0\\
  0&\cdots&0&1&0
 \end{pmatrix},
 \qquad
 S A_\delta S^{-1}=tP,\qquad
 P=
 \begin{pmatrix}
  0&0&\cdots&0&1\\
  1&0&\cdots&0&0\\
  0&1&\ddots&\vdots&\vdots\\
  \vdots&\ddots&\ddots&0&0\\
  0&\cdots&0&1&0
 \end{pmatrix}.
\]
Indeed, if \(\mathbf x'=A_\delta\mathbf x\) and
\(\mathbf y'=S\mathbf x'\), then we can compute directly
\begin{align*}
 y'_0&=x_{n-1}=t\,y_{n-1},
 y'_1=t^{1-n}\pi_F x_0=t\,y_0,\\
 y'_j&=t^{j-n}x_{j-1}
       =t^{j-n}t^{n-j+1}y_{j-1}
       =t\,y_{j-1}\qquad(2\leq j<n).
\end{align*}
Thus \(\mathbf y'=tP\mathbf y\), so projectively
\begin{equation*}
 \delta[y_0:y_1:\cdots:y_{n-1}]
   =[t y_{n-1}:t y_0:\cdots:t y_{n-2}]
   =[y_{n-1}:y_0:\cdots:y_{n-2}].
\end{equation*}
Consequently \(U_i=\delta^i\Omega_0\) and
\(\bigcup_{i=0}^{n-1}U_i=X\).  With
\(\sigma_0=(\pi_K|_\Delta)^{-1}\), define
\begin{equation*}
 \sigma_i=\delta^i\circ\sigma_0\circ\delta^{-i},
 \qquad \pi_K\circ\sigma_i=\operatorname{id}_{U_i}.
\end{equation*}
Since the action of \(D^\times\) commutes with \(G\), the same \(K\) works
for every \(i\).  The sections on \(U_I\) is induced by restriction.
\end{proof}

With \(z_{ij}=y_j/y_i\), toric coordinates on \(U_i\) are
\begin{equation}\label{eq:cover-toric-coordinates}
 T_{ij}=1+\pi_F z_{ij}\quad(j\ne i),\qquad
 U_i\simeq\{ |T_{ij}-1|\leq|\pi_F|\ (j\ne i)\}
       \subset\operatorname{Spa}
          \bigl(C\langle T_{ij}^{\pm1}\rangle,
                \cO_C\langle T_{ij}^{\pm1}\rangle\bigr).
\end{equation}
Indeed, \(|T_{ij}|=1\) and \(z_{ij}=(T_{ij}-1)/\pi_F\).
These coordinates give toric charts on \(U_I\), and the sections \(\sigma_i\) also restrict to \(U_I\).

For a continuous Banach \(E\)-representation \(\Pi\) of \(G\),
fix a \(K\)-stable open bounded lattice \(\Lambda\subset\Pi\)
(Remark~\ref{rem:Banach-no-admissibility}) and put
\(\Lambda_s:=\Lambda/p^s\Lambda\).  Let \(U\subset X\)
be a rational affinoid with a section
\begin{equation*}
 \sigma:U\longrightarrow\cM_K,
 \qquad \pi_K\circ\sigma=\operatorname{id}_U.
\end{equation*}
Rational subdomains of the \(U_i\), with restricted sections,
form a basis, as used before \cite[Lemma~3.7]{Scholze}.
Call \((U,\sigma)\) a chart, abbreviated to \(\sigma\) in
subscripts.  For \(U'\subset U\), set
\(\sigma':=\sigma|_{U'}\).

\begin{remark}
\label{rem:Banach-no-admissibility}
For a Banach norm \(\|\cdot\|\) on \(\Pi\), by joint continuity
there is an open neighborhood \(V\ni1\) and \(C_V>0\) with
\(\|h\|_{\mathrm{op}}\leq C_V\) for \(h\in V\).
Now we take the supremum over a finite cover \(K=\bigcup_{a=1}^r k_aV\), to get
\begin{equation*}
 \sup_{h\in K}\|h\|_{\mathrm{op}}
 \leq C_V\max_{1\leq a\leq r}\|k_a\|_{\mathrm{op}}
 =:C_K<\infty.
\end{equation*}
Consequently
\begin{equation*}
 \|v\|_K:=\max_{h\in K}\|hv\|,
 \qquad \|v\|\leq\|v\|_K\leq C_K\|v\|
\end{equation*}
is an equivalent non-archimedean norm, with
\(\|kv\|_K=\|v\|_K\) for \(k\in K\).  Its unit ball
\begin{equation*}
 \Lambda:=\{v\in\Pi:\|v\|_K\leq1\}
\end{equation*}
is a \(K\)-stable open bounded lattice, and by completeness we can write
\begin{equation*}
 \Lambda\simeq\varprojlim_s\Lambda/p^s\Lambda,\qquad
 \Pi=\Lambda[1/p].
\end{equation*}
The compact group $K$ acts smoothly on $\Lambda_s$, thus
for finite \(A\subset\Lambda_s\),
\begin{equation*}
 N_A=(\cO_E/p^s)\langle K\cdot A\rangle\text{ finite},
 \qquad \Lambda_s=\varinjlim_A N_A.
\qedhere
\end{equation*}
\end{remark}

\begin{lemma}
\label{lem:RH-local-torsor-formula}
For every qcqs object \(W\to U\) of \(U_{\proet}\), there is a
natural isomorphism
\begin{equation*}
 \cF_\Pi|_U(W)
 \xrightarrow{\sim}
 C_K
 \left(
  \left|W\times_UT_\sigma\right|,
  \Pi
 \right).
\end{equation*}
Here
\begin{equation*}
 T_\sigma
 \coloneqq U\times_{\sigma,\cM_K}\cM,
 \qquad
 p_\sigma:T_\sigma\longrightarrow U.
\end{equation*}
\end{lemma}

\begin{proof}
The pro-finite-\'etale right \(K\)-torsor \(T_\sigma\to U\) gives, with
$(t,g)\cdot k=(tk,k^{-1}g)$,
\[
 \Phi_\sigma:T_\sigma\times^K G
 \rightarrow U\times_X\cM,
 [t,g]\longmapsto(p_\sigma(t),tg).
\]
Over $T_\sigma$, this is the $K$-equivariant isomorphism of torsors
\[\begin{aligned}
T_\sigma\times_{U}(T_\sigma\times^KG)\cong T_\sigma\times G= U\times_{\sigma, \cM_K}\cM&\times G\cong T_\sigma\times_X\cM=T_\sigma\times_{U}(U\times_X\cM)\\
\text{via }([t_1, t_2], g)&\mapsto(t_1, t_2g).
\end{aligned}\]
Hence $\Phi_\sigma$ is an isomorphism by the faithfully flat descent along $T_\sigma\to U$.

Put $T_{\sigma,W}:=W\times_U T_\sigma$ and choose representatives
$g_c$ for $c\in K\backslash G$.  Since $K$ is open in $G$, we can identify
the geometric and topological contracted products by
\[
 |T_{\sigma,W}\times^K G|
 \cong\coprod_{c\in K\backslash G}|T_{\sigma,W}|
 \cong\bigl(|T_{\sigma,W}|\times G\bigr)/K.
\]
The maps to the $c$-summand and their inverses are given by
\[
 [t,kg_c]\longmapsto(c,tk),\qquad
 (c,t)\longmapsto[t,g_c].
\]

We define the following maps by restriction and equivariant extension:
\begin{align*}
 \operatorname{ev}_1:C_G(|T_{\sigma,W}\times^K G|,\Pi)
 &\longrightarrow C_K(|T_{\sigma,W}|,\Pi),
 &F&\longmapsto(t\mapsto F([t,1])),\\
 \operatorname{ext}:C_K(|T_{\sigma,W}|,\Pi)
 &\longrightarrow C_G(|T_{\sigma,W}\times^K G|,\Pi),
 &f&\longmapsto([t,g]\mapsto g^{-1}f(t)).
\end{align*}
These maps are continuous, inverse to each other, and natural in $W$, proving the assertion.
\end{proof}

We now show that each \(\cF_\Pi|_{U_i}\) is a complete unitary
pro-\'etale \(\Qp\)-local system in the sense of
Definition~\ref{def:complete-unitary-proetale-local-system}.

\begin{definition}
\label{def:local-arithmetic-coefficients}
\label{def:general-unitary-local-system}
For the fixed lattice $\Lambda$ and a finite $K$-stable
submodule $N\subset\Lambda_s$, we define the following sheaves on $U_{\proet}$ by their values on $W\to U$:
\begin{align*}
 \cF_{\Pi,\sigma}^+(W)
 &=C_K\bigl(|W\times_UT_\sigma|,\Lambda\bigr),\\
 \cF_{\Pi,\sigma,s}^+(W)
 &=C_K\bigl(|W\times_UT_\sigma|,\Lambda_s\bigr),\\
 \cF_{\Pi,\sigma,s,N}^+(W)
 &=C_K\bigl(|W\times_UT_\sigma|,N\bigr).
\end{align*}
\end{definition}
When \(W\to U\) is quasi-compact quasi-separated, the underlying topological space \(|W\times_UT_\sigma|\) is
quasi-compact, so continuous \(\Pi\)-valued functions are
bounded.  Thus
\begin{equation}\label{eq:local-lattice-rationalization}
 \cF_{\Pi,\sigma}^+[1/p]\simeq\cF_\Pi|_U, \qquad\cF_{\Pi,\sigma,s}^+=\varinjlim_{N\subset\Lambda_s}\cF_{\Pi,\sigma,s, N}^+.
\end{equation}
For a rational subset \(U'\subset U\) and \(\sigma'=\sigma|_{U'}\), we write
\begin{equation}\label{eq:local-coefficient-restriction}
 T_{\sigma'}=U'\times_U T_\sigma,\quad
 \cF_{\Pi,\sigma'}^+=\cF_{\Pi,\sigma}^+|_{U'},\quad
 \cF_{\Pi,\sigma',s}^+=\cF_{\Pi,\sigma,s}^+|_{U'}.
\end{equation}
The same identities hold under \'etale pullback.

\begin{proposition}
\label{lem:local-unitary}
\label{lem:general-unitary-completion}
For every chart \((U,\sigma)\), the sheaf \(\cF_\Pi|_U\) is a complete
unitary pro-\'etale \(\Qp\)-local system, with lattice
\(\cF_{\Pi,\sigma}^+\) as in
Definition~\ref{def:complete-unitary-proetale-local-system}.
\end{proposition}

\begin{proof}
Write \(\nu:U_{\proet}\to U_{\et}\).
The group \(K\) is compact, acting continuously on \(\Lambda_s\), therefore
\(\Lambda_s=\varinjlim_NN\) is a filtered colimit of finite \(K\)-stable submodules $N$.
If we denote \(K_N=\ker(K\to\operatorname{Aut}(N))\) and 
\[
 r_N:T_\sigma/K_N\longrightarrow U,\qquad
 \mathcal L_{\sigma,s,N}
 :=\bigl(r_{N,*}\underline N\bigr)^{K/K_N}.
\]
The map \(r_N\) is a finite \'etale \(K/K_N\)-torsor, so
\(\mathcal L_{\sigma,s,N}\) is finite locally constant.
For qcqs \(W\to U\), continuous maps
\(\lvert W\times_UT_\sigma\rvert\to\Lambda_s\) have finite image, giving
\begin{equation}\label{eq:general-unitary-ind-finite}
 \cF_{\Pi,\sigma,s}^+\cong\nu^*\mathcal L_{\sigma,s},\qquad
 \mathcal L_{\sigma,s}:=\varinjlim_N\mathcal L_{\sigma,s,N}.
\end{equation}
The torsion comparison \cite[Corollary~3.17]{ScholzeRigid} then implies
\(R\nu_*\cF_{\Pi,\sigma,s}^+\simeq\mathcal L_{\sigma,s}\).

On \(T_\sigma\), a set-theoretic section
\(\Lambda_s\to\Lambda\) continuously lifts
\(\Lambda_s\)-valued functions. By torsion-freeness of $\cF_{\Pi,\sigma}^+$, we obtain
\begin{gather}
 0\longrightarrow\cF_{\Pi,\sigma}^+
 \xrightarrow{p^s}\cF_{\Pi,\sigma}^+
 \longrightarrow\cF_{\Pi,\sigma,s}^+\longrightarrow0,
 \label{eq:RH-local-sheaf-reduction}\\
 0\longrightarrow\cF_{\Pi,\sigma,1}^+
 \xrightarrow{p^s}\cF_{\Pi,\sigma,s+1}^+
 \longrightarrow\cF_{\Pi,\sigma,s}^+\longrightarrow0.
 \label{eq:RH-local-successive-reduction}
\end{gather}
Since $\cF_{\Pi, \sigma}^+$ is flat over $\Zp$, we get \( \cF_{\Pi,\sigma}^+\otimes_{\Zp}^L(\Zp/p^s) \simeq\cF_{\Pi,\sigma,s}^+.\)

Let \(V\to U\) be a connected affinoid \'etale chart and
\(\widetilde V\to V\) be its universal pro-finite-\'etale cover.
By the proof of \cite[Theorem~4.9]{ScholzeRigid},
\(\widetilde V\) is affinoid perfectoid and
\[
 H^q_{\proet}(\widetilde V,\underline N)=0
 \qquad(q>0,\quad N\text{ finite }p\text{-primary}).
\]
The objects \(W=\widetilde V\times Z\), with \(Z\) profinite,
refine a pro-\'etale basis of affinoid perfectoid objects which trivialize \(T_\sigma\).
And on these basis we have 
\[
 R\Gamma_{\proet}(W,\cF_{\Pi,\sigma,s}^+)=C(Z,\Lambda_s)[0],\qquad
 C(Z,\Lambda_{s+1})\twoheadrightarrow C(Z,\Lambda_s).
\]
Computing locally and applying \cite[Lemma~3.18]{ScholzeRigid}, we get
\begin{equation}\label{eq:general-unitary-completion}
 \cF_{\Pi,\sigma}^+
 =\varprojlim_s\cF_{\Pi,\sigma,s}^+
 \simeq R\!\varprojlim_s\cF_{\Pi,\sigma,s}^+.
\end{equation}
Together with \eqref{eq:local-lattice-rationalization}, we obtain
all the conditions of Definition~\ref{def:complete-unitary-proetale-local-system}.
\end{proof}

We henceforth use \(\cF_{\Pi,\sigma,s}^+\) also for its
\'etale direct image.

\begin{definition}
\label{def:local-three-coefficients}
\label{def:local-structural-coefficients}
\label{def:local-period-coefficients}
For \(\mathscr B\in\{\Qp,\widehat{\mathscr O},\mathbb B_J\}\), put
\[
 (\mathscr B^+,t)=
 \begin{cases}
 (\Zp,p),&\mathscr B=\Qp,\\
 (\widehat{\mathscr O}^{+},p),&\mathscr B=\widehat{\mathscr O},\\
 (\mathbb A_J,a_J),&\mathscr B=\mathbb B_J.
 \end{cases}
\]
For the cases $\mathscr{B}=\widehat{\mathscr{O}}$ or $\mathbb{B}_J$, we use the almost settings specified before
Proposition~\ref{prop:geometric-almost-Cartan-Leray}. We define sheaves on \(U_{\proet}\)
\begin{align*}
 \cF_{\Pi,\sigma}^{\mathscr B}&:=\cF_\Pi^{\mathscr B}|_U,\\
 \cF_{\Pi,\sigma,s}^{\mathscr B,+}
 &:=\cF_{\Pi,\sigma}^+\otimes_{\Zp}
       (\mathscr B_U^+/t^{s}),\\
 \cF_{\Pi,\sigma}^{\mathscr B,+}
 &:=\varprojlim_s\cF_{\Pi,\sigma,s}^{\mathscr B,+}.
\end{align*}
The finite-level tensor product is the ordinary tensor
product of sheaves.  In particular,
\[
 \cF_{\Pi,\sigma,s}^{\Qp,+}=\cF_{\Pi,\sigma,s}^+,
 \qquad
 \cF_{\Pi,\sigma,s}^{\widehat{\mathscr O},+}
 =\cF_{\Pi,\sigma,s}^+\otimes_{\Zp/p^s}
       \widehat{\mathscr O}_U^+/p^s.
\]
\end{definition}

By Proposition~\ref{lem:local-unitary}, \(\cF_\Pi|_U\) is a complete
unitary pro-\'etale \(\Qp\)-local system with lattice \(\cF_{\Pi,\sigma}^+\).
We apply Lemmas~\ref{lem:structural-local-completion} and
\ref{lem:period-local-completion} to obtain
\begin{gather*}
 \cF_{\Pi,\sigma}^{\Qp,+}
 =\cF_{\Pi,\sigma}^+
 \simeq R\!\varprojlim_s\cF_{\Pi,\sigma,s}^{\Qp,+},\\
 \cF_{\Pi,\sigma}^{\mathscr B,+}
 \acong R\!\varprojlim_s
       \cF_{\Pi,\sigma,s}^{\mathscr B,+}
 \quad(\mathscr B=\widehat{\mathscr O},\mathbb B_J),\\
 \cF_{\Pi,\sigma}^{\mathscr B,+}[1/p]
 \simeq\cF_\Pi^{\mathscr B}|_U
 \quad(\mathscr B=\Qp,\widehat{\mathscr O},\mathbb B_J).
\end{gather*}
We will use the almost setting \((\cO_C,\mathfrak m_C)\) for
\(\widehat{\mathscr O}\), and
\((A_{\inf}(\cO_{C^\flat}),\mathfrak m_{\inf})\) for
\(\mathbb B_J\).  The corresponding ideal becomes
the unit ideal in \(\mathscr B\).
Consequently, for every qcqs object \(W\to U\),
\begin{equation}\label{eq:local-three-coefficient-completion}
 R\Gamma_{\proet}(W,\cF_\Pi^{\mathscr B})
 \simeq\left(R\!\varprojlim_s R\Gamma_{\proet}
       (W,\cF_{\Pi,\sigma,s}^{\mathscr B,+})\right)[1/p].
\end{equation}
All these constructions commute with inherited restriction:
\[
 \sigma'=\sigma|_{U'}\quad\Longrightarrow\quad
 \cF_{\Pi,\sigma',s}^{\mathscr B,+}
 =\cF_{\Pi,\sigma,s}^{\mathscr B,+}|_{U'},\qquad
 \cF_{\Pi,\sigma'}^{\mathscr B,+}
 =\cF_{\Pi,\sigma}^{\mathscr B,+}|_{U'}.
\]
\subsection{Cohomological amplitude for all three constructions}
\label{subsec:three-coefficient-amplitude}
\label{subsec:arithmetic-amplitude}
\label{subsec:structural-sharp-amplitude}
\label{sec:structural-period-amplitude}
\label{subsec:period-amplitude}
By Proposition~\ref{lem:local-unitary}, \(\cF_\Pi|_U\) is a complete
unitary pro-\'etale \(\Qp\)-local system with lattice \(\cF_{\Pi,\sigma}^+\).
We can therefore apply Theorem~\ref{thm:general-unitary-coefficient-amplitude}.

\begin{lemma}
\label{lem:local-three-coefficient-amplitude}
\label{lem:compact-local-Artin}
\label{lem:arithmetic-local-completed-amplitude}
For every affinoid \'etale \(V\to U\), we have
\begin{equation}\label{eq:local-geometric-almost-amplitude}
 \begin{gathered}
 R\Gamma_{\proet}(V,\cF_{\Pi,\sigma}^+)
 \in D^{[0,n-1]},\\
 R\Gamma_{\proet}(V,\cF_{\Pi,\sigma}^{\mathscr B,+})
 \ain D^{[0,n-1]}
 \qquad
 (\mathscr B=\widehat{\mathscr O},\mathbb B_J).
 \end{gathered}
\end{equation}
Here we use almost mathematics with the settings of
Definition~\ref{def:local-three-coefficients}.
\end{lemma}
\begin{proof}
      By Proposition~\ref{lem:local-unitary}, \(\cF_\Pi|_V\) is a complete
unitary pro-\'etale \(\Qp\)-local system.  The case $\mathscr{B}=\mathbb{Q}_p$ follows from
Theorem~\ref{thm:general-local-lattice-amplitude}, and the remaining cases from Proposition~\ref{prop:geometric-almost-Cartan-Leray}
with local bound \(n-1\), since \(V\) is \'etale over an affinoid
chart in \(X=\Period\).
\begin{equation}
 R\Gamma_{\proet}(V,\cF_{\Pi,\sigma}^{+})
 \in D^{[0,n-1]},\quad R\Gamma_{\proet}(V,\cF_{\Pi,\sigma}^{\mathscr B,+})
 \ain D^{[0,n-1]}
 \qquad
 (\mathscr B=\widehat{\mathscr O},\mathbb B_J).
\end{equation}
\end{proof}

We can now prove the following main result.

\begin{theorem}
\label{thm:general-Banach-amplitude}
\label{thm:structural-amplitude}
\label{thm:period-amplitude}
For every continuous Banach \(E\)-representation \(\Pi\) of
\(G\), every compact rational interval \(J\subset(0,\infty)\),
and every \(\mathscr B\in\{\Qp,\widehat{\mathscr O},\mathbb B_J\}\),
\[
 \JL_{\mathscr B}(\Pi)
 =R\Gamma_{\proet}(X,\cF_\Pi^{\mathscr B})
 \in D^{[0,2n-2]}(\Rep_E^\square(D^\times)).
\]
\end{theorem}
\begin{proof}
For each nonempty \(I\subset\{0,\ldots,n-1\}\), choose
\(i\in I\) and put \(\sigma_I=\sigma_i|_{U_I}\), using
Proposition~\ref{prop:GH-rational-cover}, and then
Lemma~\ref{lem:local-three-coefficient-amplitude} above proves the local bound on each $U_I$
\[
 R\Gamma_{\proet}(U_I,\cF_\Pi^{\mathscr B})\in D^{[0,n-1]}.
\]
The \v Cech spectral sequence is written as
\[
 E_1^{a,b}=\bigoplus_{|I|=a+1}
 H^b_{\proet}(U_I,\cF_\Pi^{\mathscr B})
 \Longrightarrow H^{a+b}\JL_{\mathscr B}(\Pi),\qquad
 E_1^{a,b}=0\quad\text{unless }0\leq a,b\leq n-1.
\]
Then the spectral sequence proves the bound; the solid assertion follows from
Theorem~\ref{thm:general-unitary-coefficient-amplitude}.
\end{proof}

\begin{remark}
\label{rem:AGV-globalization}
For every \(\mathscr B\in\{\Qp,\widehat{\mathscr O},\mathbb B_J\}\),
on each chart $U_i$ we have 
\[
 (R\nu_*\cF_\Pi^{\mathscr B})|_{(U_i)_{\et}}
 \simeq R\nu_*(\cF_\Pi^{\mathscr B}|_{U_i}),\qquad
 R^b\nu_*\cF_\Pi^{\mathscr B}=0\quad(b>n-1),
\]
by sheafifying the relation $R\Gamma_{\proet}(V,\cF_\Pi^{\mathscr B})
 \in D^{[0,n-1]}$
on the affinoid \'etale basis.  By
\cite[Corollary~2.4.13]{AyoubGallauerVezzani},
\begin{equation}\label{eq:rational-etale-cd-input}
 H^a_{\et}(X,\mathcal G)=0\qquad
 \left(a>\dim X,\quad
 \mathcal G\in\operatorname{Sh}(X_{\et},E_{\mathrm{abs}})\right)
\end{equation}for qcqs rigid variety $X$ over $C$ with Krull dimension $\dim X$ and any sheaf $\mathcal G$ of $E_{\mathrm{abs}}$-vector spaces on $X_{\et}$.
Here \(E_{\mathrm{abs}}\) is the abstract coefficient field.
The Leray spectral sequence \cite[Tag~0732]{Stacks} gives us
\[
 E_2^{a,b}=H^a_{\et}(X,R^b\nu_*\cF_\Pi^{\mathscr B})
 \Longrightarrow H^{a+b}\JL_{\mathscr B}(\Pi),\qquad
 E_2^{a,b}=0\quad\text{unless }0\leq a,b\leq n-1.
\]
This recovers Theorem~\ref{thm:general-Banach-amplitude}.
One may use the same argument to prove that for Banach or LB local systems \(\mathcal L\)
on \(X_{\proet}\) that admit an integral lattice locally for the rational topology. 
The cohomological amplitude of $R\Gamma_{\proet}(X, \cL)$ would be bounded by \(2\dim X\).
\end{remark}

\begin{remark}
\label{rem:solid-amplitude-tests}
This remark will discuss the passage to solid cohomology.
For a profinite set \(Z\), let \(\underline Z\) be its constant
pro-\'etale sheaf over \(C\).  On qcqs \(W\to X\),
\begin{equation*}
 \cF_\Pi(W)(Z)
 =C_G\bigl(|(W\times\underline Z)\times_X\cM|,\Pi\bigr).
\end{equation*}
This is \(G\)-equivariant descent of \(\underline \Pi_{\cM}\)
\cite[Lemma~4.4(i) and Corollary~4.5(i)]{AnschutzMannDescent},
\cite[\S7.2]{DRC}; use later condition \eqref{eq:LB-coefficient-continuity}
for LB representations.

On \(T_\sigma\), the sheaves $\cF_{\Pi, \sigma}$ and $\cF_{\Pi,\sigma}^+$ are trivialized by
\(\underline\Pi\) and \(\underline\Lambda\). Evaluating on an extremally
disconnected profinite \(Z\), we get \(C(Z,\Pi)\) and \(C(Z,\Lambda)\),
respectively. Since $Z$ is compact, we have
\begin{equation*}
 C(Z,\Lambda)/p^s
 \cong C(Z,\Lambda/p^s),
 \qquad
 C(Z,\Lambda)
 \cong\varprojlim_s C(Z,\Lambda/p^s),
\end{equation*}
The reductions are smooth \(K\)-modules, hence filtered unions
of finite stable submodules.  Evaluation commutes with derived
sections and \(R\nu_*\), giving the bound for the solid complex when $\mathscr{B}=\mathbb{Q}_p$.
For \(\mathscr B=\widehat{\mathscr O},\mathbb B_J\), we know from
Lemma~\ref{lem:general-geometric-perfectoid-parameters} that
\[
 R\Gamma_{\proet}(W\times Z,\cF_\Pi^{\mathscr B})
 \simeq C\!\left(Z,H^0(W,\cF_\Pi^{\mathscr B})\right)[0]
 \qquad(W\text{ affinoid perfectoid}).
\]
These evaluations prove the solid assertion of
Theorem~\ref{thm:general-unitary-coefficient-amplitude}.
For LB representations, take their filtered colimit over Banach stages locally, before \v Cech totalization.
\end{remark}

\subsection{LB and locally analytic representations}
\label{subsec:LB-lalg-amplitude}
The local cohomology bound in Lemma~\ref{lem:local-three-coefficient-amplitude}
is uniform over all Banach spaces and \(\mathscr B\).
Since filtered colimits commute with local cohomology, we can pass from
the Banach representation $\Pi$ to an LB representation under the following condition.

\noindent\textbf{The LB condition.}
Let \(\Pi\) be a Hausdorff locally convex \(E\)-vector space with a
continuous \(G\)-action.  We require a countable LB presentation
\begin{equation}\label{eq:LB-Banach-presentation}
 \Pi=\varinjlim_{m\geq1} B_m,\qquad
 B_m\hookrightarrow B_{m+1},
 \qquad B_m\text{ an }E\text{-Banach space},
\end{equation}
with injective transitions and, for every profinite \(Z\),
\begin{equation}\label{eq:LB-coefficient-continuity}
 \varinjlim_m C(Z,B_m)\xrightarrow{\sim}C(Z,\Pi).
\end{equation}
The \(B_m\) need not be \(G\)-stable; the second condition means
\(\underline\Pi\simeq\varinjlim_m\underline{B_m}\) as condensed objects.

\begin{remark}
Compact transitions imply \eqref{eq:LB-coefficient-continuity}: given a pro-finite set $Z$,
a continuous \(v:Z\to\Pi\) has compact image in a Banach
subspace \(B\), with its induced topology
\cite[Proposition~1.1.11]{EmertonAnalytic}.  By
\cite[Propositions~1.1.2(ii) and~1.1.10]{EmertonAnalytic},
\(B\to\Pi\) factors continuously through some \(B_m\).
This proves surjectivity; injectivity follows from
\(B_m\hookrightarrow\Pi\).
\end{remark}

\begin{theorem}
\label{thm:LB-amplitude}
\label{cor:structural-LB-amplitude}
\label{cor:period-LB-amplitude}
Suppose that \(\Pi\) satisfies the LB conditions
\eqref{eq:LB-Banach-presentation} and \eqref{eq:LB-coefficient-continuity}.
For every \(\mathscr B\in\{\Qp,\widehat{\mathscr O},\mathbb B_J\}\)
and every compact interval \(J\subset(0,\infty)\) with rational endpoints, we have
\begin{equation*}
 \JL_{\mathscr B}(\Pi)
 \in D^{[0,2n-2]}(\Rep_E^\square(D^\times)).
\end{equation*}
In particular, this applies to every continuous representation of compact type.
\end{theorem}

Before proving the theorem, we write each LB representation of a compact group $H$
as a colimit of Banach representations of $H$.

\begin{lemma}
\label{lem:LB-compact-Banach-stages}
Assume the LB conditions \eqref{eq:LB-Banach-presentation} and
\eqref{eq:LB-coefficient-continuity}, and fix a compact subgroup
\(H\subset G\).  There are continuous Banach representations
\(\Pi_m\) of \(H\), with injective transition maps and continuous
\(H\)-equivariant injections \(\Pi_m\to\Pi\), such that
\begin{equation}\label{eq:LB-equivariant-coefficient-continuity}
 \varinjlim_m C(Z,\Pi_m)\xrightarrow{\sim}C(Z,\Pi)
 \qquad(Z\text{ profinite}).
\end{equation}
\end{lemma}
\begin{proof}
As in \cite[Proposition~3.2.14]{EmertonAnalytic}, we write
\(i_m:B_m\hookrightarrow\Pi\) and give \(C(H,B_m)\) its
supremum norm.  Set
\begin{equation*}
 \Pi_m=\bigl\{f\in C(H,B_m):
       i_m(f(h))=h\,i_m(f(1))\text{ for every }h\in H\bigr\}.
\end{equation*}
The maps \(f\mapsto i_m(f(h))-h\,i_m(f(1))\) are continuous
with Hausdorff target, so \(\Pi_m\) is closed and Banach.  The formulas
\begin{equation*}
 (k\cdot f)(h)=f(hk),\qquad
 e_m(f)=i_m(f(1))
\end{equation*}
define an isometric continuous \(H\)-action on \(\Pi_m\) and an
injective equivariant map \(e_m:\Pi_m\to\Pi\).  

For \(v\in C(Z,\Pi)\), the map
\begin{equation*}
 Z\times H\longrightarrow\Pi,\qquad (z,h)\longmapsto h\,v(z)
\end{equation*}
is continuous, so \eqref{eq:LB-coefficient-continuity} lifts it
to \(\widetilde v\in C(Z\times H,B_m)\) for some \(m\).
The exponential law
\begin{equation*}
 C(Z\times H,B_m)\cong C\bigl(Z,C(H,B_m)\bigr)
\end{equation*}
identifies \(\widetilde v\) with a map \(Z\to \Pi_m\), proving
surjectivity of the statement.
Injectivity follows from \(\Pi_m\hookrightarrow\Pi\).
\end{proof}

\begin{proof}[Proof of Theorem~\ref{thm:LB-amplitude}]
We choose the stages \(\Pi_m\) of
Lemma~\ref{lem:LB-compact-Banach-stages} above for \(H=K\), and fix
\(K\)-stable open bounded lattices \(\Lambda_m\subset \Pi_m\).
On a chart \((U,\sigma)\), we define
\begin{align*}
 \cF_{m,\sigma}(W)
 &:=C_K(|W\times_UT_\sigma|,\Pi_m),\\
 \cF^+_{m,\sigma}(W)
 &:=C_K(|W\times_UT_\sigma|,\Lambda_m),\\
 \cF^+_{m,\sigma,s}(W)
 &:=C_K(|W\times_UT_\sigma|,\Lambda_m/p^s\Lambda_m),\\
 \cF_{m,\sigma}^{\mathscr B}
 &:=\cF_{m,\sigma}\otimes_{\Qp}^{L,\square}\mathscr B|_U,
 \qquad\mathscr B\in\{\Qp,\widehat{\mathscr O},\mathbb B_J\}.
\end{align*}
Then \(\cF_{m,\sigma}^{\Qp}=\cF_{m,\sigma}\). 
Then using the same argument of Proposition~\ref{lem:local-unitary}, 
one can easily deduce that each $\cF_{m, \sigma}$ is a complete unitary pro-\'etale $\mathbb{Q}_p$ local system 
with lattice $\cF_{m,\sigma}^+$ using the datum above. 
Thus on each $V\to U$ \'etale, we have the local bound by Lemma~\ref{thm:general-unitary-coefficient-amplitude} that
\[R\Gamma_{\proet}(V, \cF^{\mathscr{B}}_{m, \sigma})\in D^{[0,n-1]},\qquad \mathscr{B}\in\{\mathbb{Q}_p,\widehat{\mathscr{O}}, \mathbb{B}_J\}.\]

Using Lemma~\ref{lem:LB-compact-Banach-stages} and compatibility
of the solid tensor product with filtered colimits, we obtain on each $V\to U$ \'etale
\begin{equation}\label{eq:LB-all-coefficient-continuity}
 \cF_\Pi^{\mathscr B}|_U
 \simeq\varinjlim_m\cF_{m,\sigma}^{\mathscr B}\quad \Rightarrow\quad R\Gamma_{\proet}(V,\cF_{\Pi}^{\mathscr{B}})=\varinjlim_m R\Gamma_{\proet}(V,\cF_{m, \sigma}^{\mathscr{B}})\in D^{[0,n-1]}.
\end{equation}

Using the local bound above on the cover $\{U_I\}$ of
Proposition~\ref{prop:GH-rational-cover}, we get the spectral sequence
\[
 E_1^{a,b}=\bigoplus_{|I|=a+1}
 H^b_{\proet}(U_I,\cF_\Pi^{\mathscr B})
 \Longrightarrow H^{a+b}\JL_{\mathscr B}(\Pi),\qquad
 E_1^{a,b}=0\quad\text{unless }0\leq a,b\leq n-1.
\]
Thus we deduce that
\[
 \JL_{\mathscr B}(\Pi)\in D^{[0,2n-2]}
 \qquad
 \bigl(\mathscr B\in\{\Qp,\widehat{\mathscr O},\mathbb B_J\}\bigr).
\]
Finally, the bound in the solid category follows by the evaluations of
Remark~\ref{rem:solid-amplitude-tests}.
\end{proof}

\begin{corollary}
\label{cor:locally-analytic-amplitude}
Let \(\Pi\) be a locally analytic \(E\)-representation of
\(G\) whose underlying locally convex space is of compact type.
This includes all admissible locally analytic representations.
For every \(\mathscr B\in\{\Qp,\widehat{\mathscr O},\mathbb B_J\}\),
with \(J\subset(0,\infty)\) any fixed compact rational interval, we have
\begin{equation*}
 \JL_{\mathscr B}(\Pi)
 \in D^{[0,2n-2]}(\Rep_E^\square(D^\times)),
\end{equation*}
\end{corollary}
\begin{proof}
A locally analytic action is continuous
\cite[Definition~3.6.9 and the subsequent discussion]{EmertonAnalytic}.
Apply the compact-type case of Theorem~\ref{thm:LB-amplitude}.
\end{proof}

\begin{remark}
\label{rem:amplitude-endpoints}
Take \(\mathscr B=\Qp\), so \(\JL_{\mathscr B}=\JL\).  For \(\Pi\) as in Theorems~\ref{thm:general-Banach-amplitude},
\ref{thm:LB-amplitude}, or
Corollary~\ref{cor:locally-analytic-amplitude}, the component
calculation of \cite[Proposition~4.7]{Scholze} gives, as
\(E\)-vector spaces,
\begin{align*}
 \JL^0(\Pi)&=C_G(|\cM|,\Pi)\\
 &\cong C_G(F^\times,\Pi)\cong \Pi^{\mathrm{SL}_n(F)}.
\end{align*}
The \(G\)-action on \(F^\times\) is by determinant; the last
map evaluates at \(1\).  Hence
\begin{equation*}
 \Pi^{\mathrm{SL}_n(F)}=0
 \quad\Longrightarrow\quad \JL(\Pi)\in D^{[1,2n-2]}.
\end{equation*}
We have \(R^b\nu_*\cF_\Pi=0\) for \(b>n-1\) by sheafification.
In Remark~\ref{rem:AGV-globalization}, \(E_2^{n-1,n-1}\)
is the unique term of total degree \(2n-2\) and survives.  Thus
\begin{equation}\label{eq:JL-top-degree-Leray}
 \JL^{2n-2}(\Pi)\cong
 H^{n-1}_{\et}(\Period,R^{n-1}\nu_*\cF_\Pi).
\end{equation}
\end{remark}

\subsection{Compatibility of the Jacquet--Langlands functor}
In this subsection, we prove several compatibilities of all three Jacquet-Langlands functors constructed above.
\paragraph{Central character} We first prove the preservation of the central characters.

\begin{proposition}
\label{prop:JL-central-character}
Let \(\Pi\) be a Banach or LB representation of
\(G=\mathrm{GL}_n(F)\) as in Section~\ref{sec:cohomological-amplitude},
with continuous central character
\begin{equation*}
 \omega:F^\times\longrightarrow E^\times,\qquad
 (zI_n)v=\omega(z)v.
\end{equation*}
Identify \(Z(G)=F^\times=Z(D^\times)\).  Then the action
of \(z\in Z(D^\times)\) on \(\cF_\Pi\) is
\(\omega(z)\operatorname{id}\).  For
\begin{equation*}
 \mathscr B\in\{\Qp,\widehat{\mathscr O},\mathbb B_J\},
 \qquad
 \JL_{\mathscr B}(\Pi):=
 R\Gamma_{\proet}\!\left(X,
       \cF_\Pi\otimes_{\Qp}^{L,\square}\mathscr B\right),
\end{equation*}
with \(J\) a fixed compact annulus,
\begin{equation*}
 (z-\omega(z))H^j\JL_{\mathscr B}(\Pi)=0
 \qquad(z\in F^\times,\ j\in\mathbb Z).
\end{equation*}
In particular, \(\JL^j(\Pi)\) has central character \(\omega\).
\end{proposition}

\begin{proof}
Write \(z_G=zI_n\) and \(z_D=z\).
With the right \(G\)-action and left \(D^\times\)-action of
Section~\ref{sec:cohomological-amplitude}, the diagonal-centre
relation \cite[proof of Theorem~8.1(2)]{TaylorBundles} is given by 
\begin{equation*}
 z_Dx=xz_G\qquad(z\in F^\times).
\end{equation*}
The centre $Z(D^\times)$ fixes \(X=\mathbb P_C^{n-1}\). Any function 
\(f\in\cF_\Pi(W)\) is \(G\)-equivariance. Hence we have
\begin{equation*}
 (z_D\cdot f)(x)
 =f(z_D^{-1}x)
 =f(xz_G^{-1})
 =z_Gf(x)
 =\omega(z)f(x).
\end{equation*}
What's more, the centre $Z(D^\times)$ acts trivially on \(\mathscr B\), since it fixes
\(X\).  Hence the center-character relation \((z-\omega(z))H^j\JL_{\mathscr B}(\Pi)=0\) follows.
\end{proof}

\paragraph{Twisting characters} Next we prove the compatibility with twisting characters. 
\begin{proposition}
\label{prop:JL-character-twist}
Let \(\Pi\) be a Banach or LB representation of \(G\) as in
Section~\ref{sec:cohomological-amplitude}, and let
\(\eta:F^\times\to E^\times\) be a continuous character. Put
\[
 \chi=\eta\circ\det,\qquad
 \chi_D=\eta\circ\operatorname{Nrd}_{D/F},\qquad
 \chi_W=\eta\circ\operatorname{Art}_F^{-1},
\]
where \(\operatorname{Art}_F^{-1}:W_F^{\mathrm{ab}}\xrightarrow{\sim}F^\times\)
is the inverse of local reciprocity normalized by arithmetic Frobenius.

We have an isomorphism of \(D^\times\times W_F\)-equivariant
solid \(E\)-module sheaves on \(X_{\proet}\), functorial in \(\Pi\):
\[
 \cF_{\Pi\otimes_E E(\chi)}
 \simeq
 \cF_\Pi\otimes_{\underline E_X}^{\square}
 \underline{E(\chi_D\boxtimes\chi_W^{-1})}_X.
\]

Therefore for \(\mathscr B\in\{\Qp,\widehat{\mathscr O},\mathbb B_J\}\),
we consequently have equivariant equivalences
\[
 \JL_{\mathscr B}(\Pi\otimes_E E(\chi))
 \simeq
 \JL_{\mathscr B}(\Pi)\otimes_E^{L,\square}
 E(\chi_D\boxtimes\chi_W^{-1}).
\]
\end{proposition}

\begin{proof}
According to \cite[Theorem~4.4]{StrauchComponents} and
\cite[\S4, discussion preceding Proposition~4.4]{MiedaExplicitLLC},
after fixing an origin in the determinant torsor, we get a
continuous map
\[
  \kappa:|\cM|\longrightarrow\pi_0(|\cM|)
 \xrightarrow{\sim}F^\times.
\]
With our right \(G\)-action, left \(D^\times\)- and \(W_F\)-actions,
the map $\kappa$ satisfies the relation
\[
 \kappa(a\,w(x)\,g)
 =
 \operatorname{Nrd}_{D/F}(a)\,
 (\operatorname{Art}_F^{-1})(w)^{-1}\,
 \kappa(x)\det(g)
 \qquad
 (g\in G,\ a\in D^\times,\ w\in W_F).
\]

Put \(h_\eta=\eta(\kappa)^{-1}\). For \(U\to X\) pro-\'etale,
write \(U_\infty=U\times_X\cM\). We define the \(E\)-linear map
\[
 \begin{aligned}
 T_{\eta,U}\colon\cF_\Pi(U)\otimes_E E(\chi_D\boxtimes\chi_W^{-1})
 &\longrightarrow\cF_{\Pi\otimes_E E(\chi)}(U),\\
 f\otimes v&\longmapsto
 \bigl(x\longmapsto h_\eta(x)f(x)\otimes v\bigr).
 \end{aligned}
\]
The map is well-defined since for \(g\in G\), we have
\[
 \begin{aligned}
 T_{\eta,U}(f\otimes v)(xg)
 &=\chi(g)^{-1}h_\eta(x)(g^{-1}f(x))\otimes v\\
 &=g^{-1}T_{\eta,U}(f\otimes v)(x).
 \end{aligned}
\]

For section transport from
\(U\) to \((a,w)U\), we get
\[
 \begin{aligned}
 \bigl((a,w)\cdot T_{\eta,U}(f\otimes v)\bigr)(x)
 &=h_\eta((a,w)^{-1}x)f((a,w)^{-1}x)\otimes v\\
 &=T_{\eta,(a,w)U}\bigl((a,w)\cdot f\otimes
       \chi_D(a)\chi_W(w)^{-1}v\bigr)(x).
 \end{aligned}
\]
Thus \(T_\eta\) is \(D^\times\times W_F\)-equivariant.
Multiplication by \(h_\eta^{-1}\) gives its inverse.
These maps respect restrictions and profinite parameters.

The remaining assertions follows directly by tensoring with $\mathscr{B}$ and taking derived sections.
\end{proof}

\paragraph{With locally analytic vectors} We prove compatibility with taking locally analytic vectors
for Banach and LB representations satisfying $\Pi^{RG-\la}=\Pi^{G-\la}[0]$.
We follow \cite[Theorem~7.2.1]{DRC}, where the proof is local.

\begin{lemma}
\label{lem:period-proet-v-comparison}
Let \(W/C\) be affinoid perfectoid, and let
\(\mathscr B=\widehat{\mathscr O}\) or \(\mathbb B_J\), for
a fixed compact interval \(J\subset(0,\infty)\) with rational
endpoints.  Restriction from the \(v\)-site, we have
\begin{equation*}
 R\Gamma_{\proet}(W,\mathscr B)
 \xrightarrow{\sim}R\Gamma_v(W,\mathscr B)
 \simeq\mathscr B(W)[0].
\end{equation*}
This solid equivalence is natural in \(W,\mathscr B\) and
respects group actions, annulus restriction, and Frobenius.
\end{lemma}
\begin{proof}
By \cite[Lemma~4.10]{ScholzeRigid}, \(\widehat{\mathscr O}\)
is pro-\'etale acyclic on \(W\); by
\cite[Proposition~8.8]{ScholzeDiamonds}, it is also \(v\)-acyclic.
For \(\mathbb B_J\), applying Lemma~\ref{lem:period-coefficient-acyclicity}
to the trivial sheaf \(\cF=\underline{\Qp}\) proves pro-\'etale acyclicity,
and \cite[Corollary~2.2.12(2)]{DRC} proves \(v\)-acyclicity.
The comparison is the identity on \(\mathscr B(W)\).
Apply the same proof to \(W\times\underline Z\) for extremally
disconnected profinite \(Z\).  By \cite[Corollary~4.10]{BoscoDrinfeld}
for \(\mathbb B_J\), and the profinite-product formula for
\(\widehat{\mathscr O}\),
\begin{equation*}
 \mathscr B(W\times\underline Z)
 \cong C(Z,\mathscr B(W)).
\end{equation*}
Exact evaluation detects the solid equivalence; all maps are natural.
\end{proof}

\begin{proposition}
\label{prop:derived-la-JL-comparison}
Let \(\Pi\) be a Banach or LB \(E\)-representation of \(G\)
as in Subsection~\ref{subsec:LB-lalg-amplitude}.  Assume the canonical
comparison is an equivalence of solid representations:
\begin{equation}
 \Pi^{RG\text{-}\la}\simeq\Pi^{G\text{-}\la}[0].
\end{equation}
Then, for each
\(\mathscr B\in\{\Qp,\widehat{\mathscr O},\mathbb B_J\}\),
there are natural equivalences of solid
\(D^\times\)-representation complexes
\begin{equation}\label{eq:derived-la-JL-comparison}
 \JL_{\mathscr B}(\Pi)^{RD^\times\text{-}\la}
 \simeq\JL_{\mathscr B}\bigl(\Pi^{RG\text{-}\la}\bigr)
 \simeq\JL_{\mathscr B}\bigl(\Pi^{G\text{-}\la}\bigr).
\end{equation}
Here \(\JL_{\Qp}=\JL\), and \(J\) is any compact interval
with rational endpoints in \((0,\infty)\).
\end{proposition}

\begin{proof}
For compact open \(K\subset G\), by
\cite[Proposition~3.2.6(1)]{RJRCSolid} and the assumption we get
\(\Pi^{RK\text{-}\la}\simeq\Pi^{G\text{-}\la}[0]\).
We choose \(K\)-stable Banach stages \(\Pi_m\) by
Lemma~\ref{lem:LB-compact-Banach-stages} (a single stage if \(\Pi\) is Banach).
Using commutation with filtered colimits and analytic radii
\cite[Proposition~3.2.6(3)--(4) and Corollary~3.2.14(2)]{RJRCSolid},
we get, in solid vector spaces,
\[
 \Pi^{G\text{-}\la}=H^0(\Pi^{RK\text{-}\la})
 \cong\varinjlim_{m,h}\Pi_m^{K\text{-}h\text{-}\an}.
\]
Here \(h\) runs through a countable cofinal system of Banach analytic
radii.  
We follow \cite[proof of Theorem~7.2.1, Step~3]{DRC}.
For \(\mathscr B=\widehat{\mathscr O}\) or \(\mathbb B_J\),
choose a compact-torsor chart \(T_\sigma\to U\) with group
\(K\), stabilized by a compact open \(K_D\subset D^\times\).
Put
\[
 A_{\sigma,\mathscr B}=R\Gamma_{\proet}(T_\sigma,\mathscr B),\qquad
 C_{U,\mathscr B}(\Pi')=R\Gamma_{\proet}(U,\cF_{\Pi'}^{\mathscr B}).
\]
Since \(\cF_{\Pi'}|_{T_\sigma}=\underline{\Pi'}\),
we combine Lemmas~\ref{lem:structural-perfectoid-computation},
\ref{lem:structural-local-completion},
and~\ref{lem:period-coefficient-acyclicity} to compute the local values
\[
 C_{U,\mathscr B}(\Pi')\simeq R\Gamma\bigl(K,R\Gamma_{\proet}(T_\sigma, \cF_{\Pi'}^{\mathscr{B}})\bigr)\simeq
 R\Gamma\bigl(K,A_{\sigma,\mathscr B}\otimes_{\Qp}^{L,\square}\Pi'\bigr),
 \qquad \Pi'=\Pi,\ \Pi^{G\text{-}\la}.
\]
Taking filtered colimit when $\Pi'$ is LB representation.
By Lemma~\ref{lem:period-proet-v-comparison} and
\cite[Theorem~6.1.1]{DRC}, we have equivalences
\[
 A_{\sigma,\mathscr B}^{RK_D\text{-}\la}
 \xleftarrow{\sim}A_{\sigma,\mathscr B}^{R(K\times K_D)\text{-}\la}
 \xrightarrow{\sim}A_{\sigma,\mathscr B}^{RK\text{-}\la}.
\]
Since \(K_D\) acts trivially on \(\Pi\), we apply
\cite[Corollary~3.2.14(3) and Theorem~6.3.4]{RJRCSolid} to get
\begin{align*}
 \alpha_{U, \mathscr{B}}\colon C_{U,\mathscr B}(\Pi)^{RK_D\text{-}\la}
 &\simeq R\Gamma\bigl(K,
 A_{\sigma,\mathscr B}^{RK_D\text{-}\la}\otimes_{\Qp}^{L,\square}\Pi\bigr)\\
 &\simeq R\Gamma\bigl(K,
 A_{\sigma,\mathscr B}^{RK\text{-}\la}\otimes_{\Qp}^{L,\square}\Pi\bigr)\\
 &\simeq R\Gamma\bigl(K,
 A_{\sigma,\mathscr B}\otimes_{\Qp}^{L,\square}\Pi^{RK\text{-}\la}\bigr)
 \simeq C_{U,\mathscr B}(\Pi^{G\text{-}\la}).
\end{align*}
The \(\alpha_{U,\mathscr B}\) commute with restriction and
\(D^\times\)-transport and thus they glue via a \v Cech complex totalization to
\[
 \alpha_{\mathscr B}:\JL_{\mathscr B}(\Pi)^{RD^\times\text{-}\la}
 \xrightarrow{\sim}\JL_{\mathscr B}(\Pi^{G\text{-}\la}),
\]
since derived locally analytic vectors preserve finite limits.

For the case $\mathscr{B}=\mathbb{Q}_p$, we use the exact sequence of
\cite[proof of Theorem~7.2.1, Step~1]{DRC},
\[
 0\longrightarrow\Qp\longrightarrow\mathbb B_{[1,p]}
 \xrightarrow{\varphi-1}\mathbb B_{[1,1]}\longrightarrow0,
\]
Together with the period comparisons, this gives the commutative
diagram of exact triangles:
\[
 \begin{array}{ccccc}
 \JL(\Pi)^{RD^\times\text{-}\la}
 &\longrightarrow&
 \JL_{\mathbb B_{[1,p]}}(\Pi)^{RD^\times\text{-}\la}
 &\xrightarrow{\ \varphi-1\ }&
 \JL_{\mathbb B_{[1,1]}}(\Pi)^{RD^\times\text{-}\la}\\[3pt]
 {\scriptstyle\alpha_{\Qp}}\downarrow&&
 {\scriptstyle\alpha_{\mathbb B_{[1,p]}}}\downarrow_{\sim}&&
 {\scriptstyle\alpha_{\mathbb B_{[1,1]}}}\downarrow_{\sim}\\[3pt]
 \JL(\Pi^{G\text{-}\la})
 &\longrightarrow&
 \JL_{\mathbb B_{[1,p]}}(\Pi^{G\text{-}\la})
 &\xrightarrow{\ \varphi-1\ }&
 \JL_{\mathbb B_{[1,1]}}(\Pi^{G\text{-}\la}).
 \end{array}
\]
Thus \(\alpha_{\Qp}\) is an equivalence.
\end{proof}

\begin{remark}
\label{rem:Pan-local-analytic-acyclicity}
The acyclicity hypothesis does not require admissibility.
In Pan's modular-curve setting, for a nonempty \(U\) in the basis of
\cite[Theorem~4.1.7]{PanI}, put
\(B_U=\mathcal O_{K^p}(U)\) and
\(\Pi_U=\mathcal O^{\la}_{K^p}(U)=B_U^{K_p\text{-}\la}\),
where \(K_p\) is a sufficiently small compact open stabilizer.
By \cite[Proposition~4.3.15]{PanI} and
\cite[Remark~1.2.1(1) and Proposition~3.3.5]{RJRCSolid},
\[
 B_U^{RK_p\text{-}\la}\simeq \Pi_U[0],\qquad
 \Pi_U^{RK_p\text{-}\la}\simeq \Pi_U[0].
\]
For the finite-level projection
\(\pi_{K_p}:X_{K^p}\to X_{K^pK_p}\), put
\[
 V_{K_p}:=\pi_{K_p}\bigl(\pi_{\mathrm{HT}}^{-1}(U)\bigr).
\]
By \cite[\S4.3.2]{PanI}, \(V_{K_p}\) is a nonempty affinoid curve
and \(\pi_{K_p}^{-1}(V_{K_p})=\pi_{\mathrm{HT}}^{-1}(U)\).
By the Tate normalized trace, we have
\[
 \mathcal O_{K^p}(U)^{K_p}
 =\mathcal O^{\la}_{K^p}(U)^{K_p}
 \cong\mathcal O_{X_{K^pK_p}}(V_{K_p}),\qquad
 \dim_C\mathcal O_{X_{K^pK_p}}(V_{K_p})=\infty.
\]
Thus these representations have infinite-dimensional compact-open
invariants already over \(C\), and are not admissible.
\end{remark}

\paragraph{Infinitesimal character}Finally, we prove the preservation of infinitesimal characters but only for $\JL_{\widehat{\mathscr{O}}}$. 
The proof is the same as that for \cite[Theorem~7.2.2]{DRC}.
\begin{proposition}
\label{prop:nonunitary-structural-infinitesimal-character}
Let \(\Pi\) satisfy the hypotheses of
Proposition~\ref{prop:derived-la-JL-comparison}.
Suppose that \(\Pi^{G\text{-}\la}\) has
infinitesimal character
\begin{equation*}
 \chi:Z_G\longrightarrow E,\qquad
 Z_G=Z\bigl(U(\Lie G\otimes_{\Qp}E)\bigr),\qquad
 (z-\chi(z))\Pi^{G\text{-}\la}=0\quad(z\in Z_G).
\end{equation*}
Let \(Z_D=Z(U(\Lie D^\times\otimes_{\Qp}E))\), and let
\(\jmath:Z_G\xrightarrow{\sim}Z_D\) be the canonical
inner-form identification.  The structural complex
\begin{equation*}
 \JL_{\widehat{\mathscr O}}(\Pi)^{RD^\times\text{-}\la}
 \simeq\JL_{\widehat{\mathscr O}}\bigl(\Pi^{G\text{-}\la}\bigr)
\end{equation*}
is a solid locally analytic \(D^\times\)-complex with
infinitesimal character \(\chi\circ\jmath^{-1}\):
\begin{equation*}
 (\jmath(z)-\chi(z))
 H^j\JL_{\widehat{\mathscr O}}\bigl(\Pi^{G\text{-}\la}\bigr)=0
 \qquad(z\in Z_G,\ j\in\mathbb Z).
\end{equation*}
\end{proposition}

\begin{proof}
Use Proposition~\ref{prop:derived-la-JL-comparison} and put
\(\Pi_{\la}=\Pi^{G\text{-}\la}\).
We use the compact-torsor charts \((U,\sigma)\) of
Proposition~\ref{prop:derived-la-JL-comparison}, with finite-level toric coordinates.
We apply Lemmas~\ref{lem:structural-perfectoid-computation}
and~\ref{lem:structural-local-completion}, the solid Banach-stage
presentation of \(\Pi_{\la}\) from that proof, and Cartan--Leray to get
\begin{equation*}
 R\Gamma_{\proet}(U,\cF_{\Pi_{\la}}^{\widehat{\mathscr O}})
 \simeq R\Gamma\left(K,
 R\Gamma_{\proet}(T_\sigma,\widehat{\mathscr O})
                    \otimes_{\Qp}^{L,\square}\Pi_{\la}\right).
\end{equation*}
Lemma~\ref{lem:period-proet-v-comparison} identifies the sections 
\(R\Gamma_{\proet}(T_\sigma,\widehat{\mathscr O})\cong R\Gamma_v(T_\sigma,\widehat{\mathscr O})=\widehat{\mathscr{O}}(T_\sigma)\).
By \cite[Proposition~3.2.9]{RJRCSolid},
\(\Pi_{\la}^{RK\text{-}\la}\simeq \Pi_{\la}[0]\).  The projection formula and
\cite[Theorems~6.1.1 and~6.2.1]{DRC} identify the cohomology group above with
\begin{equation*}
 R\Gamma\left(K,
 \mathscr O^{\la}_{\cM}(T_\sigma)
                         \otimes_{\Qp}^{L,\square}\Pi_{\la}\right).
\end{equation*}
By \cite[Theorem~6.3.1]{DRC}, the horizontal actions agree.
By \cite[proof of Theorem~7.2.2]{DRC}, \(\jmath(z)\) acts on
this complex through the action of \(z\) on \(\Pi_{\la}\).
For \(E\)-linearity, use
\begin{equation*}
 \mathscr O^{\la}_{\cM}(T_\sigma)
       \otimes_{\Qp}^{L,\square}\Pi_{\la}
 \simeq
 \bigoplus_{\iota:E\hookrightarrow C}
 \mathscr O^{\la}_{\cM}(T_\sigma)
       \otimes_{E,\iota}^{L,\square}\Pi_{\la}.
\end{equation*}
On each summand, extend \(\chi\) along \(\iota\).  Thus
\begin{equation*}
 \jmath(z)=\chi(z)\operatorname{id}
 \quad\text{on }
 R\Gamma\left(K,
 \mathscr O^{\la}_{\cM}(T_\sigma)
                         \otimes_{\Qp}^{L,\square}\Pi_{\la}\right).
\end{equation*}
Compatibility with restriction and finite \v Cech descent
give the asserted infinitesimal character.
\end{proof}

The following example shows that primitive comparison can fail
for nonunitary representations. Compare with \cite[Example~1.1.2]{ALM}.

\begin{example}
\label{ex:standard-representation-primitive-failure}
Take \(F=E=\Qp\), \(n=2\), and let \(\Pi=\Qp^2\) with the
standard action of \(G=\mathrm{GL}_2(\Qp)\). We show that
\[
 \JL(\Pi)\otimes_{\Qp}^{L,\square}C
 \longrightarrow \JL_{\widehat{\mathscr O}}(\Pi)
\]
is not a quasi-isomorphism.

Put \(X=\Pone_C\). By \cite[Example~1.1.2]{ALM},
\(\cF_\Pi\) is the universal rational Tate module on \(X\), and
\[
 \JL^0(\Pi)=\Pi^{\mathrm{SL}_2(\Qp)}=0.
\]

We compute the underlying \(C\)-vector spaces of $\JL^0_{\widehat{\mathscr{O}}}(\Pi)$.
Let \(\nu:X_{\proet}\to X_{\et}\) be the natural morphism, and put
\[
 \widehat{\cO_X(i)}
 :=\nu^*\cO_X(i)\otimes_{\nu^*\cO_X}\widehat{\mathscr O}.
\]
We have the relative Hodge--Tate sequence \cite[Proposition~5.1.6(iii) and \S6.1]{ScholzeWeinsteinModuli},
\[
 0\longrightarrow\widehat{\cO_X(1)}(1)
 \longrightarrow\cF_\Pi\otimes_{\Qp}\widehat{\mathscr O}
 \longrightarrow\widehat{\cO_X(-1)}
 \longrightarrow0.
\]
Here the twist outside the hat is a Tate twist.
Since \(\nu_*\widehat{\cO_X(i)}=\cO_X(i)\) and
\(H^0(X,\cO_X(-1))=0\), we get
\[
 \JL_{\widehat{\mathscr O}}^0(\Pi)
 \cong H^0(X,\cO_X(1))(1)
 \cong C^2(1).
\]

Since \(C\) is solid flat over \(\Qp\)
\cite[\S2.3.1]{DRC}, the comparison map in degree zero is given by
\[
 H^0\!\left(\JL(\Pi)\otimes_{\Qp}^{L,\square}C\right)
 =0
 \longrightarrow
 \JL_{\widehat{\mathscr O}}^0(\Pi)
 \cong C^2(1).
\]
Thus primitive comparison is not an isomorphism even in
degree zero.
\end{example}

\section{Continuous principal series and the parabolic quotient}
\label{sec:continuous-parabolic-induction}
We now extend the parabolic-quotient calculation of \cite[\S5.2--5.3]{JLH}
from unitary principal series to arbitrary continuous principal series. 

Let \(B\subset G\) be the upper-triangular Borel and
\(P=P_{(n-1,1)}\supset B\).  For a continuous character
\(\eta:B\to E^\times\), put
\begin{equation*}
 \Sigma_\eta=\Ind_B^P(\eta),\qquad
 I(\eta)=\Ind_B^G(\eta)\cong\Ind_P^G(\Sigma_\eta).
\end{equation*}
Induction is continuous and unnormalized, with its Banach topology:
\begin{equation*}
 I(\eta)=\{f:G\to E\text{ continuous}:f(gb)=\eta(b)^{-1}f(g)\},
 \qquad (g_0f)(g)=f(g_0^{-1}g).
\end{equation*}
The same convention applies to \(\Ind_P^G(\Sigma_\eta)\).
Write
\begin{equation*}
 Y=\cM/P,\qquad p_P:\cM\longrightarrow Y,\qquad
 q:Y\longrightarrow X=\Period.
\end{equation*}
By \cite[Theorem~A, Proposition~4.2.4 and Corollary~A.1.2]{JLH},
\(Y\) is perfectoid and proper over \(C\), and \(|Y|\) is a spectral space
of Krull dimension \(n-1\).  In particular, \(q\) is proper and qcqs.
We define the sheaf \(\mathbb L_\eta\) on \(Y_{\proet}\) by
\begin{equation*}
 \mathbb L_\eta(W):=C_P\bigl(|W\times_Y\cM|,\Sigma_\eta\bigr),
 \qquad W\in Y_{\proet},
\end{equation*}
where \(\Sigma_\eta\) has its Banach topology.
Let \(\lambda_Y:Y_{\proet}\to |Y|\) be the morphism to the
site of open subsets of \(Y\).

We again use the method of Section~\ref{sec:cohomological-amplitude}
to compute the sheaf on rational open subsets.
Put \(P_0=P\cap K\).  Fix one \(P_0\)-stable open bounded
lattice \(\Sigma_\eta^+\subset\Sigma_\eta\), and take the
fixed \(K\)-stable lattice for \(I(\eta)\) to be
\begin{equation*}
 \Lambda_\eta=\Ind_{P_0}^{K}\Sigma_\eta^+,
 \qquad N_s=\Sigma_\eta^+/p^s\Sigma_\eta^+,
 \qquad \Lambda_s=\Lambda_\eta/p^s\Lambda_\eta.
\end{equation*}
On a chart \((U,\sigma)\),
\begin{equation*}
 G/P\cong K/P_0,\qquad Y_U=q^{-1}U\cong T_\sigma/P_0.
\end{equation*}
The map \(K\to K/P_0\) has local analytic sections.  A finite
clopen refinement gives a continuous section; lifting functions on
finite clopen partitions therefore gives us the presentation
\begin{equation*}
 \Lambda_s\cong\Ind_{P_0}^K N_s.
\end{equation*}
We give $N_s$ the discrete topology and $\Sigma_\eta^+$ its
$p$-adic topology. Thus we can define the torsion and integral sheaves.
\begin{definition}
\label{def:local-parabolic-arithmetic-coefficients}
For
\(W\in(Y_U)_{\et}\), define the \'etale sheaf
\(\mathbb L_{\eta,\sigma,s}\) by
\begin{equation*}
 \mathbb L_{\eta,\sigma,s}(W):=
 C_{P_0}\bigl(|W\times_{Y_U}T_\sigma|,N_s\bigr).
\end{equation*}
For \(W\in(Y_U)_{\proet}\), define
\begin{equation*}
 \mathbb L_{\eta,\sigma}^+(W):=
 C_{P_0}\bigl(|W\times_{Y_U}T_\sigma|,\Sigma_\eta^+\bigr).
\end{equation*}
\end{definition}
For \(\nu:(Y_U)_{\proet}\to(Y_U)_{\et}\), using the same argument as in lemma~\ref{lem:local-unitary}, one can easily get
\begin{equation*}
 \mathbb L_{\eta,\sigma}^+
 \cong\varprojlim_{s\ge1}\nu^*\mathbb L_{\eta,\sigma,s},
 \qquad
 \mathbb L_{\eta,\sigma}^+[1/p]
 \cong\mathbb L_\eta|_{Y_U}.
\end{equation*}
When taking pro-\'etale cohomology, write
\(\mathbb L_{\eta,\sigma,s}\) also for its pullback by \(\nu\).
For \(\sigma'=\sigma|_{U'}\),
\begin{equation*}
 \mathbb L_{\eta,\sigma'}^+
 =\mathbb L_{\eta,\sigma}^+|_{Y_{U'}},\qquad
 \mathbb L_{\eta,\sigma',s}
 =\mathbb L_{\eta,\sigma,s}|_{Y_{U'}}.
\end{equation*}

\subsection{The Bound for JL(\texorpdfstring{$I(\eta)$}{ })} 
\begin{lemma}
\label{lem:arithmetic-parabolic-Shapiro}
For every continuous character \(\eta:B\to E^\times\), there
are canonical \(D^\times\)-equivariant equivalences
\begin{equation*}
 Rq_{\proet,*}\mathbb L_\eta\simeq\cF_{I(\eta)},
 \qquad
 \JL(I(\eta))\simeq R\Gamma_{\proet}(Y,\mathbb L_\eta).
\end{equation*}
On each chart \((U,\sigma)\), we have the torsion comparison
\begin{equation*}
 Rq_{U,\et,*}\mathbb L_{\eta,\sigma,s}
 \simeq\nu_*\cF_{I(\eta),\sigma,s}^+.
\end{equation*}
\end{lemma}\begin{proof}
Evaluation at $1$ defines the canonical comparison
\begin{equation}\label{eq:arithmetic-parabolic-comparison-map}
 c_\eta:\cF_{I(\eta)}
 \longrightarrow Rq_{\proet,*}\mathbb L_\eta,
 \qquad F\longmapsto[t\longmapsto F(t)(1)].
\end{equation}
The same construction applies to the local torsion coefficients.

We first prove the torsion comparison. Fix a chart $(U,\sigma)$ and $s\geq1$.
Take the perfectoid basis $W=\widetilde V\times Z$ trivialization $\cF_{I(\eta)}$ 
from the proof of Proposition~\ref{lem:local-unitary}. Then
\[
 T_\sigma\times_UW\simeq W\times K,
 \qquad
 Y_W:=Y_U\times_UW\simeq W\times K/P_0.
\]
The quotient $K\twoheadrightarrow K/P_0$ has local continuous sections, 
and by a clopen decomposition of $K/P_0$, it admits a continuous section $K/P_0\to K$. 
Fix such a continuous section $\tau:K/P_0\to K$ and then the $P_0$-torsor
$W\times K\to Y_W=W\times K/P_0$ splits along this section.
Thus the affinoid perfectoid $Y_W$ trivializes the torsion sheaf
\[
 \mathbb L_{\eta,\sigma,s}|_{Y_W}\simeq\underline{N_s}.
\]
Since $Z\times K/P_0$ is profinite, $Y_W$ is an affinoid perfectoid, the same acyclicity
calculation as in Proposition~\ref{lem:local-unitary} gives
\[
 H^i_{\proet}(Y_W,\mathbb L_{\eta,\sigma,s})=0
 \qquad(i>0).
\]

In degree zero, we have the canonical identification
\begin{align*}
 H^0(Y_W,\mathbb L_{\eta,\sigma,s})
 &\cong C(Z\times K,N_s)^{P_0}\\
 &\cong C\bigl(Z,\Ind_{P_0}^K N_s\bigr)
 = C(Z,\Lambda_s)
 \cong\cF_{I(\eta),\sigma,s}^+(W),
\end{align*}
where $(h\cdot f)(z,k)=h f(z,kh)$.
The comparison sends $F$ to $(z,k)\mapsto F(z)(k)$ and
respects the $K$-action by left translation. Thus we have an isomorphism of sheaves
\[
 \cF_{I(\eta),\sigma,s}^+
 \xrightarrow{\sim}
 Rq_{U,\proet,*}\mathbb L_{\eta,\sigma,s}.
\]
By étale--pro-étale comparison
\cite[Proposition~14.8 and Corollary~16.7]{ScholzeDiamonds},
we obtain
\[
 Rq_{U,\et,*}\mathbb L_{\eta,\sigma,s}
 \simeq\nu_*\cF_{I(\eta),\sigma,s}^+.
\]

We now recover the rational comparison.
The transition maps
$C(Z,\Lambda_{s+1})\to C(Z,\Lambda_s)$ are surjective.
The completion argument of Proposition~\ref{lem:local-unitary},
using \cite[Lemma~3.18]{ScholzeRigid}, therefore applies.
Cohomology on the coherent objects $W$ and $Y_W$ commutes
with filtered colimits
\cite[Proposition~3.12]{ScholzeRigid}.
Hence
\begin{align*}
 R\Gamma_{\proet}(Y_W,\mathbb L_\eta)
 &\simeq
 \left(R\!\varprojlim_s
 R\Gamma_{\proet}(Y_W,\mathbb L_{\eta,\sigma,s})\right)[1/p]\\
 &\simeq C(Z,\Lambda_\eta)[1/p][0]
 \simeq C(Z,I(\eta))[0]\\
 &\simeq R\Gamma_{\proet}(W,\cF_{I(\eta)}).
\end{align*}
These identifications are induced by $c_\eta$.

Thus $c_\eta$ is an equivalence on this basis and hence globally.
It is $D^\times$-equivariant by construction.
Applying $R\Gamma_{\proet}(X,-)$ proves the remaining assertion.
\end{proof}

\begin{lemma}
\label{lem:perfectoid-arithmetic-degree}
Let \(W/C\) be affinoid perfectoid.  A filtered union \(\mathcal N\)
of finite locally constant \(p\)-power-torsion \'etale sheaves
satisfies
\(
 H^j_{\et}(W,\mathcal N)=0\) for \(j>1\).
\end{lemma}
\begin{proof}
By the tilting equivalence \cite[Theorem~7.12]{ScholzePerfectoid}, we can work
on \(W^\flat_{\et}\cong W_{\et}\).  For a finite \(\mathbb F_p\)-local system
\(\mathcal N\), we choose a finite \'etale cover $W'$ of $W$ such that $\mathcal{N}|_{W'}$ is trivial.
Then \((\mathcal N\otimes_{\mathbb F_p}\mathscr O_{W^\flat})|_{W'}\) is a finite free sheaf.
By finite \'etale descent, $\mathcal N\otimes_{\mathbb F_p}\mathscr O_{W^\flat}$ is a finite projective module on $W^\flat$.  
The purity theorem in \cite[Theorem~7.9 and Proposition~7.13]{ScholzePerfectoid}
then gives \(H^j_{\et}(W^\flat,\mathcal N\otimes_{\mathbb F_p}\mathscr O_{W^\flat})=0\) for \(j>0\).
Hence the twisted Artin--Schreier sequence
\begin{equation*}
 0\longrightarrow\mathcal N
 \longrightarrow\mathcal N\otimes_{\mathbb F_p}\mathscr O_{W^\flat}
 \xrightarrow{1\otimes\varphi-1}
 \mathcal N\otimes_{\mathbb F_p}\mathscr O_{W^\flat}
 \longrightarrow0
\end{equation*}
then gives us \(H^j_{\et}(W,\mathcal N)=0\) for \(j>1\).
A simple d\'evissage extends this to finite \(p\)-power torsion.
\end{proof}

We apply this bound to the sheaf $\mathbb{L}_\eta$ by taking inverse limits and inverting $p$.
\begin{lemma}
\label{lem:parabolic-arithmetic-pushforward}
The direct image $R\lambda_{Y,*}\mathbb L_\eta$ satisfies
\[
 R^b\lambda_{Y,*}\mathbb L_\eta=0\quad(b>1).
\]
\end{lemma}
\begin{proof}
Fix a chart \((U,\sigma)\), and let \(W\subset Y_U\) be an
affinoid perfectoid open.  Compactness of \(P_0\) and smoothness
of \(N_s\) express \(\mathbb L_{\eta,\sigma,s}|_W\) as a
filtered union of finite local systems, as in the proof of
Theorem~\ref{thm:general-Banach-amplitude}. Therefore by 
Lemma~\ref{lem:perfectoid-arithmetic-degree} above we have
\begin{equation*}
 A_{W,s}:=R\Gamma_{\et}(W,\mathbb L_{\eta,\sigma,s})
 \in D^{[0,1]}.
\end{equation*}
From the exact sequence, we get
\begin{equation*}
 0\longrightarrow\mathbb L_{\eta,\sigma,1}
 \xrightarrow{p^s}\mathbb L_{\eta,\sigma,s+1}
 \longrightarrow\mathbb L_{\eta,\sigma,s}\longrightarrow0
\Rightarrow
 H^1(A_{W,s+1})\longrightarrow H^1(A_{W,s})
 \longrightarrow H^2(A_{W,1})=0.
\end{equation*}
Thus \(R^1\!\varprojlim_s H^1(A_{W,s})=0\).
The Milnor sequence \eqref{eq:general-lattice-Milnor} gives us
\begin{equation*}
 H^2\!\left(R\!\varprojlim_s A_{W,s}\right)
 \cong R^1\!\varprojlim_s H^1(A_{W,s})=0,
 \qquad R\!\varprojlim_s A_{W,s}\in D^{[0,1]}.
\end{equation*}
Applying the completion argument of
Lemma~\ref{lem:local-unitary} to the \(P_0\)-torsor, together with
torsion comparison and filtered colimits on the qcqs object \(W\), we get
\begin{equation*}
 R\Gamma_{\proet}(W,\mathbb L_\eta)
 \simeq\left(R\!\varprojlim_s
 R\Gamma_{\et}(W,\mathbb L_{\eta,\sigma,s})\right)[1/p]
 \in D^{[0,1]}.
\end{equation*}
The \(Y_{U_i}\) cover \(Y\), and their affinoid perfectoid
opens form a basis of \(|Y|\).  Hence for every \(y\in|Y|\),
\begin{equation*}
 (R^b\lambda_{Y,*}\mathbb L_\eta)_y
 \cong\varinjlim_{y\in W\subset Y_{U_i}}
 H^b_{\proet}(W,\mathbb L_\eta)=0\qquad(b>1).
\end{equation*}
This proves the assertion.
\end{proof}

We can finally prove the main result.

\begin{theorem}
\label{thm:continuous-principal-series-amplitude}
For every continuous \(\eta:B\to E^\times\), we have the cohomological amplitude
\begin{equation*}
 \JL(I(\eta))\in D^{[0,n]}.
\end{equation*}
\end{theorem}
\begin{proof}
As in \cite[proof of Theorem~5.3.1]{JLH}, the spectral space
\(|Y|\) has cohomological dimension at most \(n-1\).
The spectral sequence
\begin{equation}
 E_2^{a,b}=H^a(|Y|,R^b\lambda_{Y,*}\mathbb L_\eta)
 \Longrightarrow\JL^{a+b}(I(\eta))
\end{equation}
is supported in \(0\leq a\leq n-1\), \(0\leq b\leq1\),
by Lemma~\ref{lem:parabolic-arithmetic-pushforward}.
Hence its abutment lies in \([0,n]\).
\end{proof}
\begin{remark}
If \(\eta(B)\subset\cO_E^\times\), by the torsion vanishing result
\cite[Theorem~5.3.1]{JLH}, d\'evissage, and the integral
comparison \cite[Corollary~7.1.4]{DRC}, we can deduce that
\[
 \JL(I(\eta))\in D^{[0,n-1]},
\]
with every \(\JL^i(I(\eta))\) an admissible Banach representation.
The unitary vanishing uses global primitive comparison
\cite[Proposition~5.2.4]{JLH}; for \(G=\mathrm{GL}_2(\Qp)\),
see \cite[Theorem~4.6]{Ludwig}. And the primitive comparison fails to be an isomorphism when $\eta$ is not unitary.

For an arbitrary continuous \(\eta:B\to E^\times\),
the spectral sequence above also gives us
\[
 \JL^n(I(\eta))
 \cong H^{n-1}(|Y|,R^1\lambda_{Y,*}\mathbb L_\eta).
\]
Thus the range improves to \([0,n-1]\) exactly when this
group vanishes. 
\end{remark}

\subsection{The Bounds for \texorpdfstring{$\JL_{\widehat{\mathscr{O}}}(I(\eta))$}{ } and \texorpdfstring{$\JL_{\mathbb{B}_J}(I(\eta))$}{ }}
\label{subsec:principal-series-geometric-coefficients}
We now compute the bounds for $\JL_{\widehat{\mathscr{O}}}(I(\eta))$ and $\JL_{\mathbb{B}_J}(I(\eta))$. 
Let \(J\subset(0,\infty)\) be a compact interval with rational
endpoints.  For \(\mathscr B=\widehat{\mathscr O}\) or
\(\mathbb B_J\), write \(\mathscr B_Y\) and \(\mathscr B_X\)
for the corresponding sheaves on \(Y_{\proet}\) and
\(X_{\proet}\), and set
\begin{equation*}
 \mathbb L_\eta^{\mathscr B}
 :=\mathbb L_\eta\otimes_{\Qp}^{L,\square}\mathscr B_Y,
 \qquad
 \mathbb L_{\eta,\sigma}^{\mathscr B}
 :=\mathbb L_\eta^{\mathscr B}|_{Y_U},
 \qquad Y_U=q^{-1}U.
\end{equation*}
As in Subsection~\ref{subsec:three-coefficient-amplitude}, tensoring
over \(\Qp\) retains the \(E\)-action without choosing an
embedding of \(E\) into the period ring.  Banach solid flatness
makes these sheaves to concentrate in degree zero.

\begin{definition}
\label{def:local-parabolic-geometric-coefficients}
Using Definition~\ref{def:local-parabolic-arithmetic-coefficients},
we define the following sheaves on $(Y_U)_{\proet}$
\begin{align*}
 \mathbb L_{\eta,\sigma,s}^{\widehat{\mathscr O},+}
 &:=\nu^*\mathbb L_{\eta,\sigma,s}
       \otimes_{\Zp/p^s}\widehat{\mathscr O}_{Y_U}^+/p^s,\\
 \mathbb L_{\eta,\sigma,s}^{\mathbb B_J,+}
 &:=\mathbb L_{\eta,\sigma}^+
       \otimes_{\Zp}(\mathbb A_{J,Y_U}/a_J^s),\\
 \mathbb L_{\eta,\sigma}^{\mathscr B,+}
 &:=\varprojlim_s\mathbb L_{\eta,\sigma,s}^{\mathscr B,+}.
\end{align*}
\end{definition}
Since \(p^s\mathbb A_J\subset a_J^s\mathbb A_J\), the period
model at level \(s\) depends only on \(N_s\).
For \(\sigma'=\sigma|_{U'}\),
\begin{equation*}
 \mathbb L_{\eta,\sigma',s}^{\mathscr B,+}
 =\mathbb L_{\eta,\sigma,s}^{\mathscr B,+}|_{Y_{U'}},
 \qquad
 \mathbb L_{\eta,\sigma'}^{\mathscr B,+}
 =\mathbb L_{\eta,\sigma}^{\mathscr B,+}|_{Y_{U'}}.
\end{equation*}

We first compute the sheaves locally, and prove an 
analogue of lemma~\ref{lem:parabolic-arithmetic-pushforward}.
\begin{lemma}
\label{lem:parabolic-geometric-perfectoid}
Let \(\eta:B\to E^\times\) be continuous and
\(\mathscr B\in\{\widehat{\mathscr O},\mathbb B_J\}\).
For every chart \((U,\sigma)\), the local models satisfy
\begin{equation*}
 \mathbb L_{\eta,\sigma}^{\mathscr B,+}
 \acong R\!\varprojlim_s
          \mathbb L_{\eta,\sigma,s}^{\mathscr B,+},
 \qquad
 \mathbb L_{\eta,\sigma}^{\mathscr B,+}[1/p]
 \simeq\mathbb L_{\eta,\sigma}^{\mathscr B}.
\end{equation*}
For every affinoid perfectoid object \(W\to Y_U\),
\begin{equation*}
 H^j_{\proet}(W,\mathbb L_{\eta,\sigma,s}^{\mathscr B,+})\acong0
 \quad(j>0),\qquad
 R\Gamma_{\proet}(W,\mathbb L_\eta^{\mathscr B})
 \simeq\Gamma(W,\mathbb L_\eta^{\mathscr B})[0].
\end{equation*}
We use the almost mathematics with respect to \((\cO_C,\mathfrak m_C)\) if
\(\mathscr{B}=\widehat{\mathscr O}\) and
\((A_{\inf}(\cO_{C^\flat}),\mathfrak m_{\inf})\) if
\(\mathscr{B}=\mathbb B_J\), as in Section~\ref{sec:general-local-amplitude}.
\end{lemma}
\begin{proof}
Fix an affinoid perfectoid \(W\to Y_U\).  As in the proof of
Lemma~\ref{lem:local-unitary}, with \(K\) replaced by \(P_0\),
\(\mathbb L_{\eta,\sigma,1}|_W\) is a filtered union of finite
local systems.  For $\mathscr{B}=\widehat{\mathscr{O}}$, apply
\cite[Lemma~4.12]{ScholzeRigid}; for $\mathscr{B}=\mathbb{B}_J$, use
the tilted calculation of Lemma~\ref{lem:period-integral-perfectoid}.
Using the same filtered colimits and the d\'evissage arguments, we get
\[
 H^j_{\proet}(W,\mathbb L_{\eta,\sigma,s}^{\mathscr B,+})\acong0
 \quad(j>0)\Rightarrow
 H^0_{\proet}(W,\mathbb L_{\eta,\sigma,s+1}^{\mathscr B,+})^a
 \twoheadrightarrow
 H^0_{\proet}(W,\mathbb L_{\eta,\sigma,s}^{\mathscr B,+})^a.
\]
The completion and rationalization arguments of
Lemmas~\ref{lem:structural-local-completion}
and~\ref{lem:period-local-completion}, applied to the
\(P_0\)-torsor \(T_\sigma\to Y_U\) and the lattice
\(\Sigma_\eta^+\), prove the asserted sheaf identifications.
Since both ideals in their almost setting become the unit ideal after inverting
\(p\), the surjectivity above yields
\begin{equation*}
 R\Gamma_{\proet}(W,\mathbb L_\eta^{\mathscr B})
 \simeq
 \left(R\!\varprojlim_s R\Gamma_{\proet}
      (W,\mathbb L_{\eta,\sigma,s}^{\mathscr B,+})\right)[1/p]
 \in D^{[0,0]}.
\end{equation*}
We obtain the bound in the solid category by applying the same argument
to \(W\times Z\), for extremally disconnected profinite \(Z\).
\end{proof}

We next prove an analogue of Lemma~\ref{lem:arithmetic-parabolic-Shapiro}.
\begin{proposition}
\label{prop:parabolic-geometric-Shapiro}
For every continuous character \(\eta:B\to E^\times\) and
\(\mathscr B\in\{\widehat{\mathscr O},\mathbb B_J\}\), the
canonical comparison maps are natural equivalences
\begin{equation*}
 \cF_{I(\eta)}^{\mathscr B}
 \xrightarrow{\sim}Rq_{\proet,*}\mathbb L_\eta^{\mathscr B},
 \qquad
 \JL_{\mathscr B}(I(\eta))
 \xrightarrow{\sim}R\Gamma_{\proet}(Y,\mathbb L_\eta^{\mathscr B}).
\end{equation*}
They are compatibile with the \(D^\times\)-action, 
restriction of annuli, and Frobenius on period sheaves.
\end{proposition}
\begin{proof}
The evaluation map of Lemma~\ref{lem:arithmetic-parabolic-Shapiro}
and coefficient sheaves pullback give us the natural morphism
\[
 q^{-1}\cF_{I(\eta)}\otimes_{\Qp}^{L,\square}q^{-1}\mathscr B_X
 \longrightarrow
 \mathbb L_\eta\otimes_{\Qp}^{L,\square}\mathscr B_Y.
\]
The adjunction of $q^{-1}\dashv q_*$ defines the stated comparison morphism.  Fix a chart $(U,\sigma)$
and an affinoid perfectoid pro-\'etale cover $W\to U$
trivializing $T_\sigma$; such objects form a covering basis.  
Choose a continuous section of $K\to K/P_0$, and then we can identify
\[
 Y_U\times_UW\simeq W\times K/P_0,\qquad
 I(\eta)|_K\cong C(K/P_0,\Sigma_\eta).
\]
By \cite[Corollary~4.10]{BoscoDrinfeld} for $\mathbb{B}_J$, and by
the affinoid ring of the profinite product for
$\widehat{\mathscr O}$,
\[
 \mathscr B(W\times K/P_0)\cong C(K/P_0,\mathscr B(W)).
\]
By Lemma~\ref{lem:parabolic-geometric-perfectoid} and the local calculations on $X$, we have equations
\begin{align*}
 R\Gamma_{\proet}(W,\cF_{I(\eta)}^{\mathscr B})
 &\simeq\bigl(C(K/P_0,\Sigma_\eta)\widehat\otimes_{\Qp}\mathscr B(W)\bigr)[0],\\
 R\Gamma_{\proet}(W\times K/P_0,\mathbb L_\eta^{\mathscr B})
 &\simeq\bigl(\Sigma_\eta\widehat\otimes_{\Qp}C(K/P_0,\mathscr B(W))\bigr)[0].
\end{align*}
Both $\Sigma_\eta$ and $\mathscr B(W)$ are Banach $\Qp$-spaces.
They are nuclear and solid-flat by
\cite[Corollary~A.61]{BoscoDrinfeld}, and then
\cite[Proposition~A.68]{BoscoDrinfeld} identifies their derived
solid tensor products with completed tensor products.
Thus, for any Banach $\Qp$-space $V$, we apply
\cite[Proposition~A.55(i)]{BoscoDrinfeld} to get
\[
 C(K/P_0,\Qp)\widehat\otimes_{\Qp}V
 \xrightarrow{\sim}C(K/P_0,V).
\]
By associativity, we get the natural identifications
\[
 C(K/P_0,\Sigma_\eta)\widehat\otimes_{\Qp}\mathscr B(W)
 \cong C(K/P_0,\Qp)\widehat{\otimes}_{\Qp}(\Sigma_\eta\widehat\otimes_{\Qp}\mathscr B(W))
 \cong\Sigma_\eta\widehat\otimes_{\Qp}C(K/P_0,\mathscr B(W)).
\]
Therefore \(\cF_{I(\eta)}^{\mathscr B}\cong Rq_{\proet,*}\mathbb{L}_\eta^{\mathscr{B}}.\)
The compatibilities with descent, the $D^\times$-action,
annulus restriction, and Frobenius all follow from naturality.  
By taking $R\Gamma_{\proet}(X,-)$ we get the cohomological comparison.
\end{proof}

We can now combine all the results above to get the final theorem.
\begin{theorem}
\label{thm:principal-series-geometric-amplitude}
For every continuous \(\eta:B\to E^\times\) and compact
interval \(J\subset(0,\infty)\) with rational endpoints,
\begin{equation*}
 \JL_{\widehat{\mathscr O}}(I(\eta)),\quad
 \JL_{\mathbb B_J}(I(\eta))
 \in D^{[0,n-1]}\bigl(\Rep_E^\square(D^\times)\bigr).
\end{equation*}
\end{theorem}
\begin{proof}
Affinoid perfectoid opens \(W\subset Y_{U_i}\) form an open basis
of \(|Y|\).  By Lemma~\ref{lem:parabolic-geometric-perfectoid},
we can compute near every \(y\in|Y|\),
\begin{equation*}
 (R^b\lambda_{Y,*}\mathbb L_\eta^{\mathscr B})_y
 =\varinjlim_{y\in W}
 H^b_{\proet}(W,\mathbb L_\eta^{\mathscr B})=0
 \qquad(b>0).
\end{equation*}
Thus \(R\lambda_{Y,*}\mathbb L_\eta^{\mathscr B}
\simeq\lambda_{Y,*}\mathbb L_\eta^{\mathscr B}\).  By
Proposition~\ref{prop:parabolic-geometric-Shapiro}, only the $b=0$ row of the Leray
spectral sequence survives
\begin{equation*}
 E_2^{a,b}=H^a\bigl(|Y|,R^b\lambda_{Y,*}
                                  \mathbb L_\eta^{\mathscr B}\bigr)
 \Longrightarrow\JL_{\mathscr B}^{a+b}(I(\eta)).
\end{equation*}
Now the space $|Y|$ is a spectral space of Krull dimension $n-1$, the spectral sequence then only has column $0\le a\le n-1$, $b=0$, thus
\begin{equation*}
 \JL_{\mathscr B}^a(I(\eta))
 \cong H^a\bigl(|Y|,\lambda_{Y,*}\mathbb L_\eta^{\mathscr B}\bigr)
 =0\qquad(a>n-1).
\end{equation*}
For the bound in solid category, we apply the same calculation with an
extremally disconnected profinite parameter as in Remark~\ref{rem:solid-amplitude-tests}.
\end{proof}

\begin{remark}
By \cite[proof of Theorem~7.2.1, Step~1]{DRC}, the
fundamental exact sequence, together with
Proposition~\ref{prop:parabolic-geometric-Shapiro}, gives us
another description of the extra degree $\JL^n(I(\eta))$: 
\begin{equation*}
 \JL(I(\eta))\longrightarrow
 \JL_{\mathbb B_{[1,p]}}(I(\eta))
 \xrightarrow{\varphi-1}
 \JL_{\mathbb B_{[1,1]}}(I(\eta))
 \longrightarrow\JL(I(\eta))[1].
\end{equation*}
The preceding bound \([0,n-1]\) for $\JL_{\mathbb{B}_J}(I(\eta))$ recovers the 
range \([0,n]\) for $\JL(I(\eta))$, with
\begin{equation*}
 \JL^n(I(\eta))\cong
 \coker\!\left(
 \JL_{\mathbb B_{[1,p]}}^{n-1}(I(\eta))
 \xrightarrow{\varphi-1}
 \JL_{\mathbb B_{[1,1]}}^{n-1}(I(\eta))\right).
\end{equation*}
This is an isomorphism of solid sheaves. The displayed cokernel is
the remaining obstruction to the \([0,n-1]\) range.
\end{remark}

\printbibliography[heading=bibintoc]
\end{document}